\documentclass[reqno, 12pt]{amsart}
\usepackage{amsfonts}
\usepackage{amssymb,amsmath}
\usepackage{amsthm}
\usepackage{upref}
\usepackage{bbm}

\makeatletter
\def\LaTeX{\leavevmode L\raise.42ex
    \hbox{\kern-.3em\size{\sf@size}{0pt}\selectfont A}\kern-.15em\TeX}
 \makeatother

\numberwithin{equation}{section}

\newtheorem{lemma}{Lemma}[section]
\newtheorem{theorem}[lemma]{Theorem} 
\newtheorem{corollary}[lemma]{Corollary}
\newtheorem{proposition}[lemma]{Proposition}

\theoremstyle{definition}

\newtheorem{definition}[lemma]{Definition}

\newtheorem{assumption}[lemma]{Assumption}

\theoremstyle{remark}

\newtheorem{remark}[lemma]{Remark}

\renewcommand{\det}{\operatorname{Det}}

 \newcommand{\supp}{\operatorname{supp}}
  \newcommand{\e}{\eqref}

\newcommand{\q}{\quad}

\newcommand{\ti}{\tilde}
\newcommand{\wt}{\widetilde}

\newcommand{\la}{\langle}
\newcommand{\ra}{\rangle}

\newcommand{\ov}{\overline}
 \renewcommand{\d}{\delta}

   \newcommand{\sgn}{\operatorname{sgn}}

  \newcommand{\rank}{\operatorname{rank}}
  \newcommand{\spect}{\operatorname{spect}}

   \newcommand{\ess}{\operatorname{ess}}

\renewcommand\Im{\operatorname{Im}}
\renewcommand\Re{\operatorname{Re}}

\newenvironment{pf}{\begin{proof}}{\end{proof}}

\def\qqq{\mathrel{\subset\mkern-15mu\lower.38ex\hbox{${\scriptscriptstyle\rightarrow}$}}}

\let\goth\mathfrak
\let\cal\mathcal

\let\Bbb\mathbb

  \DeclareMathOperator{\spec}{spec}
  
\begin{document}

 \title {Quasi-Carleman  operators and their spectral properties}
 \author{ D. R. Yafaev}
\address{ IRMAR, Universit\'{e} de Rennes I\\ Campus de
  Beaulieu, 35042 Rennes Cedex, FRANCE}
\email{yafaev@univ-rennes1.fr}
\keywords{ The sigma-function,  positivity,    necessary and sufficient conditions, total numbers of negative eigenvalues,  the Carleman operator,  singular Hankel  operators, the Laplace transform.}
\subjclass[2000]{47A40, 47B25}

%  \date{\today}

\begin{abstract}
The Carleman operator is defined as integral operator with kernel $(t+s)^{-1}$ in the space $L^2 ({\Bbb R}_{+}) $. This is the simplest example of a Hankel operator which can be explicitly diagonalized. Here we study a class of self-adjoint  Hankel operators   (we call them quasi-Carleman  operators) generalizing the  Carleman operator   in various directions. We   find explicit formulas for 
the total number of negative eigenvalues of quasi-Carleman  operators and, in particular, necessary and sufficient conditions for their positivity. Our approach relies on the concepts of the sigma-function and of the quasi-diagonalization of Hankel operators introduced in the preceding paper of the author.
  \end{abstract}

\maketitle

% \thispagestyle{empty}

%************************************************************
\section{Introduction}  
%***********************************************************

{\bf 1.1.}
Hankel operators can be defined as integral operators  
\begin{equation}
(H f)(t) = \int_{0}^\infty h(t+s) f(s)ds 
\label{eq:H1}\end{equation}
in the space $L^2 ({\Bbb R}_{+}) $ with kernels $h$ that depend  on the sum of variables only. We refer to the books \cite{NK,Pe,Po} for basic information on Hankel operators. Of course $H$ is symmetric if $  h(t)=\ov{h(t)}$. There are very few cases when Hankel operators can be explicitly diagonalized. The most simple and important case $h(t)=t^{-1}$ was considered by T.~Carleman in \cite{Ca}.

Here we study a class of Hankel operators (quasi-Carleman operators) with  kernels  
\begin{equation}
h(t)= (t+r)^{-q} e^{-\alpha t},  \q  r\geq 0,
\label{eq:E1r}\end{equation}
where $\alpha$ and   $q$ are real numbers.  We will see that a Hankel operator with kernel \e{eq:E1r} can be correctly defined as a self-adjoint operator in the Hilbert space $L^2 ({\Bbb R}_{+}) $ if
\begin{equation}
{\rm either}\q \alpha>0\q {\rm or}\q \alpha=0, \: q>0,
\label{eq:aq}\end{equation}
that is, $h(t)\to 0$ as $t\to \infty$. The singularity of $h(t) $ at the point $t=0$ may be arbitrary.
There are of course no chances to explicitly find the spectrum and eigenfunctions of quasi-Carleman operators. The only exceptions are the cases  $q=1$, $\alpha=0$   and $q=1$, $r=0$ considered by F.~G.~Mehler \cite{Me} and W.~Magnus \cite{Ma}, respectively (see also  \S 3.14 of the book \cite{BE} and the papers \cite{Ro}, \cite{Y1}).

%Hankel operators $H$ with such kernels are not in general bounded, but under our assumptions they are semibounded from below.  
 
To obtain information about spectral properties of quasi-Carleman operators, we  here  use
 the method   of quasi-diagonalization of Hankel operators suggested in \cite{Yqd}. Roughly speaking, this method relies on the identity
\begin{equation} 
H= {\sf L}^* \Sigma {\sf L} 
 \label{eq:MAIDp}\end{equation} 
 where ${\sf L}$ is
the Laplace transform 
defined by the relation
   \begin{equation}
  ({\sf L} f) (\lambda)= \int_{0}^\infty e^{-t \lambda} f(t) dt  
\label{eq:LAPj}\end{equation}
and $\Sigma$ is the multiplication  operator by a function $\sigma (\lambda)$.  We use the term ``sigma-function" of the Hankel operator $H$ (or of the kernel $h(t)$) for this function. It is {\it formally}  linked to the kernel $h(t)$ of $H$ by the relation
  \begin{equation}
  h (t)  =    \int_{-\infty}^\infty     e^{-t \lambda} \sigma(\lambda)  d \lambda,
   \label{eq:conv1}\end{equation}
   that is, $h (t)$ is the two-sided Laplace transform of $\sigma(\lambda)$. We consider here ${\sf L}$ as a mapping of appropriate spaces of test functions so that ${\sf L}^*$ is the corresponding mapping of dual spaces (of distributions).

It is clear from   formula \e{eq:conv1} that $\sigma(\lambda)$  can be a regular function only for   kernels $h(t)$   satisfying some specific analytic  conditions. Without such very restrictive assumptions,  $\sigma$ is necessarily a distribution. For example, for  kernels \e{eq:E1r} the sigma-function    is given by the explicit formula
    \begin{equation}
\sigma (\lambda)=
   \frac{1}{\Gamma(q)}   (\lambda -\alpha)_{+}^{q-1}e^{-r (\lambda-\alpha)} , \q -q\not\in{\Bbb Z}_{+},
\label{eq:bbr7}\end{equation}
where $\Gamma(\cdot)$ is the gamma function. The function  $\mu_{+}^{q-1}$ is   regular  for $q>0$ (in this case $\mu_{+}^{q-1}=\mu^{q-1}$ for $\mu>0$ and $\mu_{+}^{q-1}=0$ for $\mu\leq 0$), but it is a singular distribution  for $q\leq 0$. For $-n-1< q < -n$ where   $n\in {\Bbb Z}_{+}$ and    a test function $\varphi (\lambda)$, it is
defined   by the standard formula 
 \begin{equation}
\int_{-\infty}^\infty  (\lambda -\alpha)_{+}^{q-1} \varphi(\lambda)d\lambda=\int_\alpha^\infty  (\lambda -\alpha) ^{q-1} \big(\varphi(\lambda)-\sum_{p=0}^n \frac{1}{p!}\varphi^{(p)}(\alpha)(\lambda -\alpha)^{p} \big)d\lambda.
\label{eq:di}\end{equation}
If $-q\in{\Bbb Z}_{+}$, then   $\sigma(\lambda)$ is a linear combination of derivatives $\d^{(p)}(\lambda-\alpha)$ of delta-functions for $p=0,1, \ldots, -q$ (note that the corresponding Hankel operator $H$ has finite rank if $\alpha>0$). Thus,  in general,    
 $\sigma (\lambda)$ is   a singular distribution so that $\Sigma$ need not   even be defined as an operator.  Therefore, instead of operators, we   work with the corresponding quadratic forms which is both more general and more convenient. So, to be precise, instead of \e{eq:MAIDp} we consider the identity
  \begin{equation}
 (H f, f) =  ( \Sigma {\sf L} f, {\sf L}f)
\label{eq:MAIDs}\end{equation} 
  on   the set   of  elements $f\in C_{0}^\infty ({\Bbb R}_{+})$ and assume only that $h\in C_{0}^\infty ({\Bbb R}_{+})'$.   
   Note that ${\sf L}$ acts as an isomorphism of 
$C_{0}^\infty ({\Bbb R}_{+})$ onto a set (denoted $\cal Y$) of analytic functions and for kernels \e{eq:E1r} quadratic forms in \e{eq:MAIDs}  are well defined  for all values of $\alpha, q\in {\Bbb R}$ and $r\geq 0$.

It follows from \e{eq:MAIDs} that the total numbers of positive $N_{+}(H)$ and negative $N_{-}(H)$ eigenvalues of the operator $H$ equal the same quantities for $\Sigma$:
   \begin{equation}
N_{\pm}(H)= N_{\pm}(\Sigma).
\label{eq:MAIDT}\end{equation} 
In particular,    $\pm H\geq 0$ if and only if $\pm \Sigma\geq 0$.  In general, we have to speak about quadratic forms 
$ (H f, f) $ for $f\in C_{0}^\infty ({\Bbb R}_{+})$  and $  ( \Sigma w, w)$ for $w\in \cal Y$ instead of the operators $H$ and $\Sigma$. Under the only assumption $h\in C_{0}^\infty ({\Bbb R}_{+})'$ the number $N_{\pm}(H)$ is defined as the maximal dimension of linear sets in $ C_{0}^\infty ({\Bbb R}_{+})$ where $ \pm (H f, f) >0$ for all $f\neq 0$.

This general construction was applied to kernels  \e{eq:E1r}  in \cite{Yqd}  where the numbers $N_{\pm}(H)$ were explicitly calculated (see Theorem~\ref{HKL} below) as a function of $q$. In particular, it turns out that $N_{\pm}(H)$ do not depend on $\alpha$ and $r$.

\medskip

{\bf 1.2.}
This paper can be considered as a continuation of \cite{Yqd}. It has two goals. The first is to define  Hankel operators  $H$ with kernels \e{eq:E1r} as  self-adjoint operators in the space $L^2({\Bbb R}_{+})$. We note that such operators are bounded in the following two cases: 

$1^0$ if $\alpha >0$, then either $r>0$  and $q$ is arbitrary  or $r=0$ and $ q \leq 1$ , 

$2^0$ if $\alpha = 0$, then either $r>0$  and $  q \geq 1$  or $r=0$ and $ q =1$.

 It is easy to see that under assumption  \e{eq:aq} all unbounded quasi-Carleman operators have kernels
    \begin{equation}
h(t)=t^{-q} e^{-\alpha t}, \q \alpha> 0,\q q>1,
\label{eq:QCs}\end{equation}
or
     \begin{equation}
h(t)=(t+r)^{-q}, \q    q >0,
\label{eq:QC}\end{equation}
where either $r>0$, $q<1$ or $r=0$, $q\neq 1$.

A study of unbounded integral operators goes back to T.~Carleman \cite{Ca} (see also Appendix~I to the book \cite{AhGl}). In particular, his general results apply to Hankel operators with kernels satisfying the condition
 \begin{equation}
\int_{t}^\infty | h(s)|^2 d s<\infty, \q \forall t>0 .
\label{eq:Carl}\end{equation} 
 This condition allows one to define $H$ as a symmetric but not as a self-adjoint operator. We note that for kernels \e{eq:QC}   condition \e{eq:Carl} is not satisfied if $q\leq 1/2$; in this case the corresponding operator $H$ is not defined   even on the set $C_{0}^\infty({\Bbb R}_{+})$.

We proceed from the identity \e{eq:MAIDs} and define $H$ in terms of the corresponding quadratic form. In view of formula \e{eq:bbr7} for kernels \e{eq:E1r}
   the identity   \e{eq:MAIDs}  reads   as
 \begin{equation}
 \int_0^\infty \int_0^\infty   \overline{f(t)} f(s) (t+s+r)^{-q } e^{-\alpha (t+s)}dtds =
  \frac{e^{\alpha r}}{\Gamma( q)} \int_{-\infty }^\infty (\lambda-\alpha)_{+}^{q-1 }e^{-r \lambda} |({\sf L}f)(\lambda) |^2 d \lambda
\label{eq:C1}\end{equation}
where $f\in C_{0}^\infty ({\Bbb R}_{+})$ is arbitrary. The following result defines $H$ as a self-adjoint operator.

\begin{theorem}\label{QFHH}
 Let the function $h(t)$ be given by formula \e{eq:E1r} where either $\alpha>0$, $q\geq 1$ or $\alpha=0$, $q >0$ $($the parameter $r\geq 0$ is arbitrary$)$. Then the form \e{eq:C1} defined on the set of functions $f\in C_{0}^\infty ({\Bbb R}_{+})$ admits the closure  in the space $L^2 ({\Bbb R}_{+})$, and it is closed on the set of  all $f \in   L^2  $ such that the integral in the right-hand side of   \e{eq:C1}  is finite.  The corresponding Hankel operator is strictly positive.  
  \end{theorem}
  
  In particular, this result applies to kernels \e{eq:QCs} and \e{eq:QC}.
  
Actually, we consider a more general problem of defining  a   Hankel operator by means of its sigma-function
given by some   measure $dM(\lambda)$ on $[0,\infty)$. 
Thus we assume that
\begin{equation}
h(t)=\int_{0}^\infty e^{-t\lambda} dM(\lambda) . 
\label{eq:Bern1}\end{equation}
Recall that, as shown by H.~Widom  in \cite{Widom}, the Hankel operator $H$ with kernel \e{eq:Bern1} is bounded if and only if the measure in \e{eq:Bern1} satisfies the condition 
\begin{equation}
M( [0,\lambda))= O(\lambda) \;   {\rm as} \; \lambda\to 0   \;   \mbox{ and as} \; \lambda\to \infty . 
\label{eq:Wid}\end{equation}
In particular,  for the Lebesgue measure $d M(\lambda)=d\lambda$ on ${\Bbb R}_{+}$, we have $h(t)=t^{-1}$  and $H$ is the Carleman operator.

 Our goal is to study the singular case  when   condition \e{eq:Wid} is not satisfied, but  $H$ can be defined as an
  unbounded positive operator via its quadratic form $(H f,f)$. We find  sufficient (and practically necessary) conditions on the measure $d M(\lambda)$ guaranteeing that the form  $(H f,f)$ defined on $C_{0}^\infty ({\Bbb R}_{+})$ admits the closure in $L^2 ({\Bbb R}_{+})$ and describe the domain of its closure. These results are deduced from properties of the Laplace transform $\sf L$ considered as  a mapping of $L^2 ({\Bbb R}_{+})$ into $L^2 ({\Bbb R}_{+}; dM(\lambda))$.
  Perhaps, the results on the Laplace transform are of independent interest.

  For quasi-Carleman operators with homogeneous kernels, we  prove the following spectral result.

 \begin{theorem}\label{Car}
Let $H$ be the   Hankel operator  with kernel $h(t)= t^{-q}$ where $q>0$ and $q\neq 1$. 
The spectrum of the operator $H$ coincides with the positive half-line, it has a constant multiplicity  and is absolutely continuous.
 \end{theorem}

    \medskip
 
 {\bf 1.3.}
 Another goal of the paper is to study perturbations of singular Hankel operators $H_{0}$ constructed in 
 Theorem~\ref{QFHH} by bounded quasi-Carleman operators $V$ with kernels
  \begin{equation}
v(t)= v_{0}(t+\rho)^k e^{-\beta t}, \q v_{0}\in {\Bbb R}, \, \beta\geq 0, \, \rho\geq 0, \,  k\in {\Bbb R} .
\label{eq:E1v}\end{equation} 
The cases $k<0$ and $k\geq 0$ turn to to be qualitatively different. According to formula \e{eq:bbr7} in the first case the sigma-function $\sigma (\lambda)=\sigma_{0} (\lambda)+\sigma_{v} (\lambda)$ of the operator $H=H_{0}+V$ belongs to the set $L^1_{\rm loc} ({\Bbb R}_{+})$. It implies that $H\geq 0$ if  $\sigma (\lambda)\geq 0$ and $H$ has infinite negative spectrum if $\sigma(\lambda)<0$ on a set of positive measure. The precise result is stated in Theorem~\ref{HKC}.

 In the case $k\geq 0$ we have the following result.

    \begin{theorem}\label{FDHC}
Let $H_{0}$ be the Hankel operator with  kernel $h_{0}(t)$   given by formula  \e{eq:E1r}   where either $q\geq 1$  for  $\alpha >0$ or $q>0$ for $\alpha=0$. Let $V$ be the  Hankel operator with kernel \e{eq:E1v}
 where    $\beta>0$ and $k \geq 0$.   Then 
    \begin{equation}
  N_{-} (H_{0} + V)=N_{-} (  V).
\label{eq:FDb}\end{equation} 
 \end{theorem}

  Thus we obtain the striking result: the total multiplicity of the negative spectrum of the operator $H = H_{0} + V$ does not depend on the operator $H_{0}$.   The inequality $N_{-} (H)\leq N_{-} ( V)$ is of course obvious because $H_{0}\geq 0$. On the contrary, the opposite inequality $N_{-} (H)\geq N_{-} (  V)$ looks surprising because the  operator $H_{0}$  may be ``much stronger" than $V$; for example, the Hankel operator with kernel $h_{0} (t)=t^{-q}$ is never compact and is unbounded unless $q=1$. Nevertheless its adding to $V$ does not change
    the total number of   negative   eigenvalues of the operator $V$. A heuristic explanation of this phenomenon can be given in terms of the sigma-functions. The sigma-function $\sigma_{0}(\lambda)$ of the  operator  $H_{0}$ is  continuous and positive  while the sign-function $\sigma_{v}(\lambda)$  of $V$ has a strong negative singularity at the point $\lambda = \beta$. 
       Therefore the sigma-functions of $H$ and $  V$ have the same negative singularity. 
      Very loosely speaking,  the supports of the  functions $\sigma_{0}(\lambda)$  and $\sigma_{v}(\lambda)$  are essentially disjoint so  that the operators $H_{0}$ and $V$ ``live in orthogonal subspaces", and hence the positive operator $H_{0}$ does not affect the negative spectrum of $  V$. 
      
      Relation \e{eq:FDb}  is also true for perturbations of singular Hankel operators by  finite rank Hankel operators. The kernels of these operators are linear combinations of functions
       \e{eq:E1v} where $k\in{\Bbb Z}_{+}$, $r=0$ and $\Re\beta>0$. 
            The sigma-function of  such kernel  \e{eq:E1v}  consists of  the delta-function and its derivatives supported at   the point   $\lambda=\beta$. If $\Im\beta\neq 0$, this sigma-function  is more
  singular than functions  \e{eq:bbr7} which impedes the proof of relation \e{eq:FDb}.

    Let us compare  the results on the negative spectrum of Hankel and differential operators. Let   ${\sf H }=D^2+ {\sf V}(x)$ be the Schr\"odinger operator in the space $L^2 ({\Bbb R})$. Suppose that ${\sf V}(x)\leq 0$.  If ${\sf V}(x)$ decays sufficiently rapidly as $|x|\to \infty$, then  $N_{-}({\sf H })<\infty$ and $N_{-}({\sf H })=\infty$ in the opposite case. Contrary to the Schr\"odinger case,   the negative spectrum of Hankel operators $H=H_{0}+V$ is  determined not by the behavior  of $v(t)$ at singular points $t=0$ and $t=\infty$ but exclusively by the corresponding sigma-functions.

     \medskip
 
 {\bf 1.4.}
 Let us briefly describe the structure of the paper. We collect necessary results of \cite{Yqd} 
   in Section~2.    Singular Hankel operators are studied in Section~3. In particular, Theorems~\ref{QFHH} and \ref{Car} are proven there. Perturbations of singular Hankel operators by bounded quasi-Carleman operators are considered in 
   Section~4. Similar results for finite rank  perturbations  are discussed in Section~5.  In particular, the results of  Sections~4 and 5  imply Theorem~\ref{FDHC}.
     Finally, in Appendix we study   the Fourier transform sandwiched by    functions one of which is unbounded. This problem is adjacent to that considered in Section~3.
   
%   The results of the type of Theorems~\ref{FDHC} and \ref{FDFR} are obtained in Sections~4 and 5, respectively.
   
  Let us introduce some standard
 notation:    ${\cal S}={\cal S} ({\Bbb R}) $ is the Schwartz  space, $\Phi$ is  the Fourier transform, 
\[
(\Phi u) (\xi)=  (2\pi)^{-1/2} \int_{-\infty}^\infty u(x) e^{ -i x \xi} dx.
\]
For spaces of test functions, for example ${\cal S}$ and  $C_{0}^\infty ({\Bbb R}_{+})$,   we denote by ${\cal S}'$ and  $C_{0}^\infty ({\Bbb R}_{+})'$ the dual classes of distributions (continuous antilinear functionals). We use the notation   ${\pmb\la} \cdot, \cdot {\pmb\ra}$ and $\la \cdot, \cdot\ra$ for   the  
  duality symbols in $L^2 ({\Bbb R}_{+})$ and $L^2 ({\Bbb R})$, respectively. They are   linear in the first argument and antilinear in the second argument.

   We often use the same notation for a function and for the operator of multiplication by this function. The Dirac function is standardly denoted $\d(\cdot)$; $\d_{n,m}$ is the Kronecker symbol, i.e., $\d_{n,n}=1$  and $\d_{n,m}= 0$  if $n\neq m$. 
   The letter $C$ (sometimes with indices) denotes various positive constants whose precise values are inessential.

  We use the special notation $\bf C$ for the Hankel operator $H$ with kernel $h(t)=t^{-1}$, that is, for the Carleman operator. Recall that $\bf C$ is bounded and it has the absolutely continuous spectrum $[0,\pi]$ of multiplicity $2$. Its sigma-function $\sigma (\lambda)$ equals $1$ for $\lambda\geq 0$ and it equals $0$ for $\lambda < 0$.

   %$\mathbbm{1}_{X}$ is the characteristic function of the set $X$.  

%************************************************************
\section{The sigma-function and  the main identity}  
%*

Here we collect some necessary results of \cite{Yqd}. In particular, we give the precise definition of the sigma-function $\sigma(\lambda)$ and discuss the main  identity \e{eq:MAIDs}.  

\medskip

{\bf 2.1.}
We work on test functions $f  \in C_{0}^\infty ({\Bbb R}_{+})$ and   require that $h $ belong to the dual space $ C_{0}^\infty ({\Bbb R}_{+})'$. Let
the set ${\cal Y}$ consist of entire functions  $\varphi(\lambda)$ satisfying,   for all $\lambda \in \Bbb C$,   bounds 
  \begin{equation}
  | \varphi (\lambda)| \leq C_{n}  (1+| \lambda |)^{-n} e^{r_{\pm} |\Re \lambda |}, \q \pm\Re\lambda\geq 0,
 \label{eq:YYZ}\end{equation}
  for all $n$ and some $r_{+}=r_{+}(\varphi)<0$; the number $r_{-}=r_{-}(\varphi)$  may be arbitrary. Thus functions in $\cal Y$ exponentially decay as $\Re \lambda \to +\infty$,  and they are exponentially bounded as $\Re \lambda \to -\infty$.
  The space $\cal Y$ is of course invariant with respect to the complex conjugation   $\varphi(\lambda)\mapsto \varphi^* (\lambda): = \ov{\varphi (\bar{\lambda})}$. By definition, $\varphi_{k} (\lambda)\to 0$ as $k\to \infty$  in $\cal Y$ if all functions $\varphi_{k} (\lambda)$ satisfy bounds \e{eq:YYZ} with the same constants $r_+ <0$, $r_{-}$, $C_{n}$ and $\varphi_{k} (\lambda)\to 0$ as $k\to \infty$  uniformly on   compact subsets of $\Bbb C$. 
  
Let the Laplace  transform  ${\sf L}$ be defined by formula \e{eq:LAPj}.
By one of the versions of the Paley-Wiener theorem (see, e.g., the book \cite{GUECH} for similar assertions),
  ${\sf L}: C_{0}^\infty ({\Bbb R}_{+}) \to {\cal Y}$ is the one-to-one continuous mapping of $ C_{0}^\infty ({\Bbb R}_{+}) $ onto $ {\cal Y}$ and the inverse mapping 
   ${\sf L}^{-1}: {\cal Y} \to  C_{0}^\infty ({\Bbb R}_{+}) $ is also continuous. In such cases we say that $\sf L$  is an isomorphism.
    Passing to the dual spaces, we see that the mapping
      $
      {\sf L}^*: {\cal Y}' \to C_{0}^\infty    ({\Bbb R}_{+})'
$
 is also an isomorphism. 
 
 % We emphasize that we write ${\sf L} ^*$ here because this operator acts in the spaces of distributions.

   Let us construct the sigma-function. 
   
    \begin{definition}\label{sigmay}
      Assume that 
  $ h \in C_{0}^\infty ({\Bbb R}_{+})'$.
   The  distribution $\sigma \in {\cal Y}'  $ defined by the formula
 \begin{equation}
\sigma=({\sf L}^*)^{-1}h  
 \label{eq:LAPL1g}\end{equation}
is called the sigma-function of the kernel $h$ or of the corresponding Hankel operator $H$.
 \end{definition}
 
According to this definition  for all  $F\in C_{0}^\infty    ({\Bbb R}_{+})$, we have   the identity
 \begin{equation}
 { \pmb\la} h , F {\pmb \ra} =   { \pmb\la}  {\sf L}^* \sigma ,   F {\pmb \ra} =   { \la} \sigma ,  {\sf L} F {  \ra}.
  \label{eq:MAIDG}\end{equation}

By virtue of \e{eq:LAPL1g},
  the kernel  $h(t)$ can be recovered
from  its sigma-function $\sigma(\lambda)$ by the formula $h={\sf L}^*\sigma$ which gives the precise sense to  formal relation \e{eq:conv1}. Thus 
   there is the one-to-one correspondence between kernels $  h\in C_{0}^\infty    ({\Bbb R}_{+})' $ and their sigma-functions $\sigma \in {\cal Y} ' $.

     Let us introduce the   Laplace convolution
\[
( \bar{f}_{1}\star f_{2})(t)=
\int_{0}^t    \overline{f_{1}(s)} f_{2}(t-s) ds  
 \] 
 of   functions $ \bar{f}_{1},  f_{2} \in C_{0}^\infty    ({\Bbb R}_{+})$. Then it formally follows from \e{eq:H1}  that
\[
(Hf_{1} ,f_{2})=  {\pmb \la} h, \bar{f}_{1}\star f_{2} {\pmb \ra} 
 \]
 where we write ${\pmb \la} \cdot, \cdot {\pmb \ra}$ instead of $( \cdot, \cdot )$ because $h$ may be a distribution.  The following result was established in \cite{Yqd}.

 \begin{theorem}\label{1}
Let   $h \in C_{0}^\infty ({\Bbb R}_{+})' $, and let $\sigma \in {\cal Y}'$ be defined by formula
  \e{eq:LAPL1g}.
  Then       the identity
 \begin{equation}
 { \pmb\la} h ,  \bar{f}_{1}\star f_{2} {\pmb \ra} =   { \la} \sigma , (  {\sf L} f_{1} )^*{\sf L} f_{2} { \ra}
  \label{eq:MAID}\end{equation}
  holds for arbitrary
   $f_{1}, f_{2}\in C_{0}^\infty    ({\Bbb R}_{+})$.
 \end{theorem}

    The identity  \e{eq:MAID} attributes a precise meaning to   \e{eq:MAIDp} or \e{eq:MAIDs}. 
    
 \medskip
 
  {\bf 2.2.}
  Let ${\Bbb H}^p_{r}$ be the Hardy space of functions analytic in the right half-plane. By the Paley-Wiener theorem the operator $(2\pi)^{-1/2}{\sf L} : L^2 ({\Bbb R}_{+})\to {\Bbb H}^2_{r}$ is unitary.  If a Hankel operator $H$ is bounded, then necessarily the sigma-function $\sigma$ belongs to the   space $({\Bbb H}^1_r)'$ dual to ${\Bbb H}^1_{r}$, and the identity  \e{eq:MAID} extends to all $f_{1}, f_{2}\in L^2 ({\Bbb R}_{+})$.

If $\supp \sigma$ belongs to the right half-plane, it  is sometimes convenient to make the exponential change of variables  and to  define the function (we call it the sign-function) 
  \begin{equation}
 s(x)=\sigma (e^{-x}), \q -\pi/2 < \Im x < \pi /2.
  \label{eq:sign}\end{equation}
  In particular, for sigma-function \e{eq:bbr7} we have
 \[
  s (x)=
   \frac{e^{\alpha r}}{\Gamma(-k)}   (e^{-x} -\alpha)_{+}^{-k-1}e^{-r e^{-x}} ;
\]
obviously $s\in {\cal S}'$. It follows from \e{eq:sign} that
  \begin{equation}
 s[u,u]:=\la s, u^* u\ra=  { \pmb\la} \sigma,  w^* w {\pmb \ra}=:\sigma [w,w]  
  \label{eq:signX}\end{equation}
if
    \begin{equation}
  u(x)=  e^{-x/2 } w(e^{-x} ).
  \label{eq:sign1}\end{equation}
   Note that $w(\lambda)$ is analytic in  the half-plane $\Re \lambda>0$ if and only if the corresponding function $u(x)$ is analytic in the strip $-\pi/2 <\Im x<\pi/2$. Moreover, the conditions
    $w\in {\Bbb H}^2_{r}$ and
    \[
\sup_{-\pi/2 <a <\pi/2}    \int_{-\infty}^\infty | u(x+ia)|^2 dx <\infty
    \]
    are equivalent.
    
    According to identity \e{eq:signX} we can work either with sigma-functions $\sigma (\lambda)$ and test functions $w(\lambda)$  or with sign-functions $s(x)$ and  test functions  $u(x)$. Both points of view are equivalent, and we frequently pass from one   to another at our convenience.
    
    To recover $f(t)$ from $w(\lambda)= ({\sf L}f)(\lambda)$, we have to invert the Laplace transform $\sf L$. To that end, we use its factorization. Let $\Gamma(z)$ be the gamma function and
     \begin{equation}
  ({\pmb \Gamma}   g)(\xi)=\Gamma(1/2+i\xi) g(\xi), \q \xi\in {\Bbb R}.
  \label{eq:LTM5g}\end{equation}
  Note that
      \begin{equation}
   |\Gamma (1/2 +i\xi| = \sqrt{\frac{\pi} {\cosh (\pi \xi)}}=: v(\xi).
\label{eq:LMMv}\end{equation}
  Put  $ (U w) (x)= e^{ x/2} w(e^x)$. Obviously the operator $U: L^2 ({\Bbb R}_{+})\to L^2 ({\Bbb R} )$ is unitary, and hence the Mellin transform ${\sf M}=\Phi U$  is also unitary. Let   ${\cal J}$, $({\cal J} g)(\xi)=g(-\xi)$, be the reflection. Then the Laplace transform factorizes as 
    \begin{equation}
{\sf L} =   {\sf M}^{-1} {\cal J} {\pmb \Gamma}  {\sf M}
 \label{eq:LAPL1}\end{equation}
so that ${\sf L}^{-1} =   {\sf M}^{-1}   {\pmb \Gamma}^{-1}  {\cal J} {\sf M}$ and hence
 \begin{equation}
 ({\sf M} f )(\xi)= \Gamma(1/2+i\xi)^{-1}({\sf M} w )(-\xi)= \Gamma(1/2+i\xi)^{-1}(\Phi  u  )(\xi). 
  \label{eq:sign2}\end{equation} 
  This formula allows us to recover $f(t)$ if either $w(\lambda)$ or $u(x)= e^{-x/2} w (e^{-x})$ are given.

 \medskip
 
  {\bf 2.3.}
 Suppose  now that $h(t)=\overline{h(t)}$ for all $t>0$, or to be more precise $\ov{{\pmb \la} h, F{\pmb \ra}} = {\pmb \la} h, \overline{F}{\pmb \ra}$ for all $F\in C_{0}^\infty    ({\Bbb R}_{+})$. Then it follows from  \e{eq:MAIDG} that  the sigma-function is also real, that is,  $\ov{{  \la} \sigma, w{  \ra}} ={  \la} \sigma, w^*{  \ra}$ for all $w\in  {\cal Y} $.

 Below we use the following natural definition.  
 
  \begin{definition}\label{hss}
     Let ${\sf h}[\varphi , \varphi ]$ be   a real quadratic form defined on  a linear set ${\sf D} $.  We denote by $N_{\pm}({\sf h}) = N_{\pm}({\sf h}; {\sf D})$ the maximal dimension of linear sets ${\cal M}_{\pm}\subset {\sf D}$   such that $\pm {\sf h} [\varphi,\varphi] > 0$   for all $\varphi\in {\cal  M}_{\pm}$, $\varphi\neq 0$.
          \end{definition}
          
     Definition~\ref{hss} means that there exists a linear set ${\cal M}_{\pm} \subset {\sf D} $, $\dim {\cal M}_{\pm}= N_{\pm}({\sf h}; {\sf D}) $,
       such that $\pm {\sf h}[\varphi , \varphi ]  > 0$   for all $\varphi\in {\cal M}_{\pm}$, $\varphi \neq 0$, and for every  linear set ${\cal M}_{\pm}' \subset {\sf D}$ with  $\dim {\cal M}_{\pm}'> N_{\pm}({\sf h}; {\sf D})$ there exists $\varphi\in {\cal M}_{\pm}'$, $\varphi \neq 0$,  such that $\pm {\sf h}[\varphi , \varphi ]  \leq 0$.

  Of course, if the set ${\sf D} $ is dense in a Hilbert space $\cal H$ and ${\sf h}[\varphi , \varphi ]$ is semibounded and closed on ${\sf D} $, then for the self-adjoint operator ${\sf  H}$ corresponding to ${\sf h}$, we have  $ N_{\pm}({\sf  H})=N_{\pm}({\sf h}; {\sf D}) $. In particular, this is true for bounded operators ${\sf  H}$.
  
 We apply Definition~\ref{hss} to the forms $h[f,f]={\pmb\la} h, \bar{f} \star f {\pmb\ra}$ on $f \in C_{0}^\infty    ({\Bbb R}_{+})$ and $\sigma[w,w]=\la \sigma, w^* w \ra$  on $w\in   {\cal Y}$.

Since ${\sf L} : C_{0}^\infty    ({\Bbb R}_{+}) \to    {\cal Y}$ is an isomorphism,
   the following assertion is a direct consequence of Theorem~\ref{1}.
 
  \begin{theorem}\label{HBx}
  Let $h \in   C_{0}^\infty    ({\Bbb R}_{+})'$. Then $\sigma = ({\sf L}^*)^{-1}h \in  {\cal Y}' $ and
  \[
  N_{\pm}(h; C_{0}^\infty    ({\Bbb R}_{+}))= N_{\pm}(\sigma; {\cal Y} ).
  \]
   In particular,     the form $\pm {\pmb\la}h,\bar{f} \star f {\pmb\ra} \geq 0$ for all
 $f\in C_{0}^\infty    ({\Bbb R}_{+})$  if and only if the form  $\pm {\la}\sigma , w^* w  {\ra} \geq 0$ for all $w \in {\cal Y} $. 
     \end{theorem}
     
Thus a Hankel operator $H$ is positive (or negative) if and only if its sigma-function $\sigma (\lambda)$ is positive (or negative).

    \medskip
 
 {\bf 2.4.}
  The following assertion (see \cite{Y}) is very convenient for calculation of the numbers $N_{\pm}(\sigma; {\cal Y}  )$. Its proof   relies on formula \e{eq:LAPL1}.

  \begin{lemma}\label{HBxe}
Suppose   that distribution \e{eq:sign} belongs to the class $  {\cal S}' $. 
  Then      
  \[
  N_{\pm}(\sigma; {\cal Y}  )= N_{\pm}(\sigma; C_{0}^\infty ({\Bbb R}_{+}) ).
  \] 
  \end{lemma}

Putting together Theorem~\ref{HBx} and Lemma~\ref{HBxe} we obtain  the following result.

    \begin{theorem}\label{HByz}
 Let $h \in C_{0}^\infty ({\Bbb R}_{+})'$.  Suppose   that the corresponding  distribution \e{eq:sign} belongs to the class $  {\cal S}' $.  Then 
    \[
   N_{\pm}(h; C_{0}^\infty ({\Bbb R}_{+}) ) =N_{\pm}(\sigma; C_{0}^\infty ({\Bbb R}_{+})).
\]
   \end{theorem}

 In applications to quasi-Carleman operators, the following
  result was   obtained also in \cite{Yqd}.

    \begin{theorem}\label{HKL}
Let $h(t)$ be given by formula \e{eq:E1r} where   $\alpha \in {\Bbb R}$ and $r\geq 0$. If $q> 0$, then  ${\pmb \la }h, \bar{f}\star f{\pmb \ra }\geq 0$ for all $f\in C_{0}^\infty ({\Bbb R}_{+} )$. If
  $q<0$ but $|q| \not \in {\Bbb Z}_{+}$, then  $N_+ (h; C_{0}^\infty ({\Bbb R}_{+}) ) =  [|q| ] /2+1$, $N_- (h; C_{0}^\infty ({\Bbb R}_{+}) ) = \infty$   for even $[|q| ]$ and  $N_- (h; C_{0}^\infty ({\Bbb R}_{+}) ) = ([|q| ]+ 1)/2$, $N_+ (h; C_{0}^\infty ({\Bbb R}_{+}) ) = \infty$ for odd $[|q| ]$. 
 \end{theorem}

This result is of course deduced from formula \e{eq:bbr7} for the sigma-function. If $q> 0$, then $\sigma(\lambda)\geq 0$. If $q< 0$, then $\Gamma (q)<0$ for $[|q|]$ even and $\Gamma (q) >0$ for $[|q|]$ odd. Therefore, for example, for even $[|q|]$, the sigma-function is continuous and negative everywhere except the point $\lambda=\alpha$ which ensures that $N_- (h; C_{0}^\infty ({\Bbb R}_{+}) )  = \infty$. The singularity at the point $\lambda=\alpha$  produces a finite number of positive eigenvalues.

% Observe that the numbers $N_{\pm} (H)$ do not depend on  $\alpha \in {\Bbb R}$ and that $\sigma   \in  {\cal S}'_{ -1/2}$ if $\alpha>0$. Therefore, by the proof of Theorem~\ref{HKL}, relation \e{eq:abc} can be used.

 %************************************************************
  %************************************************************
\section{Singular quasi-Carleman operators}  
 %************************************************************
 %****

   {\bf 3.1.}
Let us now  consider    Hankel operator with kernels \e{eq:Bern1} where
the (locally finite nonnegative) measure $dM(\lambda)$ on $[0,\infty)$ satisfies the condition
   \begin{equation}
 \int_{0}^\infty e^{-t\lambda} d M(\lambda)<\infty,\q \forall t >0.
\label{eq:Bern1x}\end{equation}
 Then  the  Hankel quadratic form admits the representation
  \begin{equation}
{\pmb \la}h, \bar{f}\star f {\pmb \ra} =\int_{0}^\infty |({\sf L}f)(\lambda)|^2 d M(\lambda), \q \forall  f \in C_{0}^\infty({\Bbb R}_{+}),
\label{eq:QF}\end{equation}
where ${\sf L}$ is the Laplace transform  defined by formula \e{eq:LAPj}.
According to formula \e{eq:bbr7} kernels \e{eq:QCs} and \e{eq:QC} satisfy assumption \e{eq:Bern1x} with $d M(\lambda)=\sigma(\lambda)d\lambda$.

For the study of the form \e{eq:QF}, we consider  ${\sf L}$ as the mapping
of $   L^2 ({\Bbb R}_{+})=: L^2 $ into $ L^2 ({\Bbb R}_{+};dM) =: L^2 ( M)$; the scalar product in the space $L^2 (M)$ will be denoted $(\cdot,\cdot)_{M}$.
Note
 that for an arbitrary $f \in L^2  $, the integral $({\sf L}f)(\lambda)$ converges   and the function $({\sf L}f)(\lambda)$  is continuous for all $\lambda>0$; moreover, 
 $$|({\sf L}f)(\lambda)| \leq (2\lambda)^{-1/2} \| f\| .$$ 
  Let $ \cal E$ (resp. ${\cal E}_{M}$) consist of functions $f\in L^2  $ (resp. $g   \in L^2 (M) $) compactly supported in ${\Bbb R}_{+}$. If $f\in \cal E$, then $({\sf L}f)(\lambda)$ is a continuous function for all $\lambda\geq 0$ and  $({\sf L}f)(\lambda)= O(e^{-c \lambda})$ with some $c=c(f)>0$ as $\lambda\to\infty$. In particular,  ${\sf L}f\in L^2 (M)$ for such $f$.

Let us also introduce the operator ${\sf L}_{*}$ formally adjoint to $\sf L$ by the equality 
 \begin{equation}
 ( {\sf L}_{*} g)(t)=  \int_{0}^\infty   e^{-t\lambda}  g(\lambda)   d M(\lambda).
\label{eq:LTM}\end{equation} 
For an arbitrary $g \in L^2 (M) $,   integral \e{eq:LTM} converges for all $t>0$,  the function $({\sf L}_{*} g)(t)$  is continuous and is bounded as $t\to \infty$.  If $g\in \cal E_{M}$, then
$({\sf L}_{*} g)(t)$ is a continuous function for all $t\geq 0$ and  $({\sf L}_{*}g)(t)= O(e^{-c t})$ with some $c=c(g)>0$ as $ t \to\infty$; in particular,  ${\sf L}_{*}g\in L^2  $.

 The following assertion is a direct consequence of the Fubini theorem.

 \begin{lemma}\label{LTMF}
 Let assumption \e{eq:Bern1x} be satisfied.
    Suppose that either $ f \in {\cal E} $ and   $g\in L^2 (M)$ or $ f \in L^2 $ and   $g\in {\cal E}_{M}$. Then
      \begin{equation}
 ({\sf L} f, g)_{M} =   (f, {\sf L}_{*} g).
\label{eq:LTM1}\end{equation} 
    \end{lemma}

Define now the operator $A_{0}:L^2  \to   L^2 (M)$    by the equality $A_{0}f= {\sf L}f$ on the domain  ${\cal D}(A_{0})=\cal E$. Let us construct its adjoint $A_{0}^* : L^2 (M) \to L^2 $. Let the set ${\cal D}_{*}\subset   L^2 (M)$ consist of   $g \in   L^2 (M)$ such that ${\sf L}_{*}  g\in L^2$. As we have seen, $\cal E_{M}\subset {\cal D}_{*}$. 
  
     \begin{lemma}\label{LTM}
     Under assumption \e{eq:Bern1x} the operator $A_{0}^*$ is given by the equality $A_{0}^* g = {\sf L}_{*} g$ on the domain ${\cal D}(A_{0}^*) ={\cal D}_{*}$. 
    \end{lemma}
    
     \begin{pf}
     It follows from  identity \e{eq:LTM1} 
     for   $ f \in {\cal D}(A_{0})$ and   $g \in   {\cal D}_{*}$ that   $   {\cal D}_{*}\subset {\cal D}(A_{0}^*)  $  and   $A_{0}^* g = {\sf L}_{*}g$ for $g \in {\cal D}_{*}$.
         Conversely, if $g \in {\cal D}(A_{0}^*) $, then $| ({\sf L} f, g)_{M} | \leq C \| f\| $ for all $ f \in {\cal D}(A_{0})$. In view again of 
     \e{eq:LTM1}, this estimate implies that $| ( f, {\sf L}_{*} g) | \leq C \| f\| $ and hence
         ${\sf L}_{*} g\in L^2$. Thus   ${\cal D}(A_{0}^*)   \subset {\cal D}_{*}$.
          \end{pf}
  
  \begin{corollary}\label{LTM1}
The operator $A_{0} $   admits the closure   if and only if $  M(\{0\} )=0$.  
  \end{corollary}
    
      \begin{pf}
      If $  M(\{0\} )=0$, then the set $\cal E_{M}$ is dense in $L^2 (M)$. Since ${\cal E}_{M}\subset {\cal D} (A_{0}^*)$, it follows that the operator $A_{0}^*$ is densely defined, or equivalently, that the operator $A_{0} $   admits the closure.     Conversely, suppose that $  M(\{0\} ) > 0$ and denote by $dM_{0} (\lambda)$ the restriction of $dM  (\lambda)$ on ${\Bbb R}_{+}$. Since $\cal E_{M}$ is dense in $L^2 (M_{0})$, we see that      $({\sf L}_{*}  g  )(t)\to 0$ as $t\to \infty$ for all $g\in L^2 (M_{0})$ and hence $({\sf L}_{*}  g  )(t)\to M(\{0\} ) g(0)$ as $t\to \infty$ for all $g\in L^2 (M)$.  It follows that
 $g(0)=0$  for all    $g\in{\cal D}_{*} $ so that   $   {\cal D}_{*}  $ is not dense in $L^2 (M)$. 
        \end{pf}

        \begin{remark}\label{QFR1}
     The choice of the set $\cal E$ in the definition of the operator $A_{0}$ is not  essential because its restriction, for example,   on  the set $C_{0} ^\infty ({\Bbb R}_{+} )$ has the same  adjoint as $A_{0}$ defined on $\cal E$.
 \end{remark}
        
        Next, we construct the second adjoint $ A_{0}^{**}$. Let the operator $A$ be defined   by the equality $A f={\sf L} f$ on the domain   ${\cal D} (A)$ which  consists of   $f \in   L^2  $ such that ${\sf L} f\in L^2(M)$. Obviously, $A_{0}\subset A$.
        
        \begin{lemma}\label{LTM2}
     Under the assumptions $  M(\{0\} )=0$ and \e{eq:Bern1x}, we have   $A_{0}^{**}\subset A $. 
    \end{lemma}
    
     \begin{pf}
     If $f \in {\cal D} (A_{0}^{**})$, then $|(f, {\sf L}_{*} g)|\leq C \| g \|_{M}$  for all $g\in {\cal D}_{*}$ and, in particular, for $g\in {\cal E}_{M}$. Using  the identity
  \e{eq:LTM1}  
for  $ f \in L^2 $ and  $g\in \cal E_{M}$, we see that $|({\sf L}f,    g)_{M}|\leq C \| g \|_{M}$ and hence
${\sf L} f\in L^2 (M)$. Thus $f \in {\cal D} (A)$ and
 $A_{0}^{**}f= {\sf L} f$ according again to  \e{eq:LTM1}. 
          \end{pf}

  \medskip 
  
   {\bf 3.2.}
The proof of the opposite inclusion $ A \subset A_{0}^{**}$ is essentially more difficult. Now we have to check relation   \e{eq:LTM1} for all $ f \in L^2 $ such that ${\sf L} f \in L^2 (M)$ and all  $g \in L^2 (M)$   such that   ${\sf L}_{*} g \in L^2 $. To that end,  we
        require that, for some $k>0$, the measure $dM(\lambda)$ satisfies the condition
           \begin{equation}
\int_{0}^\infty (1 + \lambda)^{-k} dM(\lambda)<\infty 
\label{eq:LTM2}\end{equation} 
which is stronger than \e{eq:Bern1x}.
Let $\chi_{n}$ be the operator of multiplication by the function $\chi_{n} (\lambda) =e^{-\ln^2\lambda/n^2}$. Consider
$ ({\sf L} f, \chi_{n}g)_{M} $. Since ${\sf L} f\in L^2 (M)$, $g\in L^2 (M)$  and $ \chi_{n} \to I$ strongly in this space, we see that
\begin{equation}
\lim_{n\to \infty} ({\sf L} f, \chi_{n}g)_{M} =   ({\sf L} f,  g)_{M} .
\label{eq:XL}\end{equation}

 Let us now show that 
 \begin{equation}
\lim_{n\to \infty} ({\sf L} f, \chi_{n}g)_{M} =   (  f, {\sf L}_{*} g) .
\label{eq:LTM3}\end{equation}
 Astonishingly, this turns out to be  a substantial problem. We put
   \begin{equation}
 \hat{\chi}_{n}=    \Phi U \chi_{n} U^{-1}\Phi^{-1}.
\label{eq:LTM6}\end{equation}
Since the operator $ U \chi_{n} U^{-1} $ acts as the multiplication by the function $e^{-x^2/n^2}$, we have
 \begin{equation}
( \hat{\chi}_{n} g)(\xi)=    \int_{-\infty}^\infty \hat{\chi}_{n} (\xi -\eta)  g(\eta) d\eta
\label{eq:LTM7}\end{equation}
where
 \begin{equation}
\hat{\chi}_{n} (\xi)= (2 \sqrt{\pi})^{-1}n    e^{-n^2 \xi ^2/4}   .
 \label{eq:LTM7x}\end{equation}

 \begin{lemma}\label{LTM3}
 Let the  operator ${\pmb\Gamma}$  be given by   formula
   \e{eq:LTM5g}.
 Then the operators 
 \begin{equation}
T_{n}= {\pmb\Gamma}^* \hat{\chi}_{n}  ({\pmb\Gamma}^*)^{-1} .
\label{eq:LTM8}\end{equation}
defined on $C_{0}^\infty ({\Bbb R})$ extend to
  bounded operators  on the space $L^2 ({\Bbb R})$.   Moreover,   their norms are bounded uniformly in $n$ and   $T_{n} \to I$  strongly as $n\to\infty$.
    \end{lemma}
    
     \begin{pf}
 Let the function $v(\xi)$ be defined by formula \e{eq:LMMv}.
     It follows from  \e{eq:LTM7}, \e{eq:LTM8}  that 
      \[
|( T_{n} g)(\xi)|\leq   v(\xi) \int_{-\infty}^\infty \hat{\chi}_{n}  (\xi -\eta)  v(\eta)^{-1}  |g(\eta)| d\eta
\]
and hence
 \begin{equation}
|( T_{n} g)(\xi)| ^2 \leq  q_{n} (\xi)   \int_{-\infty}^\infty \hat{\chi}_{n}  (\xi -\eta)    |g(\eta)|^2 d\eta
\label{eq:LMM1}\end{equation}
where
 \[
q_{n} (\xi)=  v(\xi)^2  \int_{-\infty}^\infty \hat{\chi}_{n}  (\xi -\eta)  v(\eta)^{-2}  d\eta.
\]

It is easy to see that
 \begin{equation}
q_{n} (\xi)\leq Q
\label{eq:LX2}\end{equation}
where the constant $Q$ does not depend on $\xi$ and $n$. Indeed,
in view of formulas \e{eq:LTM7x} and \e{eq:LMMv}, we have  to estimate the expression 
 \begin{equation}
n v(\xi)^2  \int_{-\infty}^\infty e^{-n^2 (\xi -\eta)^2 /4  }e^{\pi\eta}  d\eta=
n v(\xi)^2 e^{\pi\xi} e^{\pi^2/n^2 } \int_{-\infty}^\infty e^{-(n  x -2 \pi /n )^2/4 }  dx
\label{eq:LX3}\end{equation}
where we   have set $ x= \eta-\xi$. Since the last integral equals $2 \sqrt{\pi}/n$, expression \e{eq:LX3} is uniformly bounded which proves \e{eq:LX2}.

Integrating now estimate \e{eq:LMM1} over $\xi$, we obtain that $\| T_{n} \| ^2 \leq Q$. 
Using that the operators $\hat{\chi}_{n} \to I$ strongly as $n\to \infty$, we see that 
$T_{n} g\to g$
 for all $g \in C_{0}^\infty ({\Bbb R})$ and hence for all $g \in L^2 ({\Bbb R})$. 
          \end{pf}
          
          Let us return to  relation \e{eq:LTM3}. 
Now we use the factorization \e{eq:LAPL1} of the Laplace operator.

                  \begin{lemma}\label{LTM3ss}
       Let the operator $T_{n} $ be defined by formula \e{eq:LTM8}, and let ${\sf M}=\Phi U$ be the Mellin transform. Then
  \begin{equation}
\chi_{n} {\sf L}  = {\sf L} {\sf M}^{-1} T_{n}^* {\sf M} .
\label{eq:LTM9}\end{equation}
    \end{lemma}
    
     \begin{pf}
    Putting together relations
\e{eq:LAPL1} (where we use that $\Phi {\cal J}= {\cal J} \Phi $),  \e{eq:LTM6} and \e{eq:LTM8}, we see that
\begin{multline}
 {\sf L} {\sf M}^{-1} T_{n}^* {\sf M} =(U^{-1} {\cal J} \Phi^{-1}  {\pmb\Gamma } {\sf M}) {\sf M}^{-1} ( {\pmb\Gamma }^{-1} \hat{\chi}_{n}  {\pmb\Gamma }){\sf M}
 \\
 =U^{-1}  {\cal J}   \Phi^{-1}  \hat{\chi}_{n} {\pmb\Gamma }{\sf M}
 =U^{-1}  {\cal J}  ( U \chi_{n} U^{-1} ) \Phi^{-1} {\pmb\Gamma }{\sf M} .
 \label{eq:LMM3L}\end{multline}
 Since the operators ${\cal J}$ and $ U \chi_{n} U^{-1} $ commute,  the right-hand side of \e{eq:LMM3L} equals
 \[
   \chi_{n} U^{-1}  {\cal J} \Phi^{-1}  {\pmb\Gamma }{\sf M} =\chi_{n} {\sf L} ,
   \]
   which proves \e{eq:LTM9}.
  \end{pf}

          It follows from identity \e{eq:LTM9} that  
\begin{equation}
({\sf L} f, \chi_{n}g)_{M} = ({\sf L} u_{n} , g)_{M}     
\label{eq:LMM2}\end{equation}
where
\begin{equation}
  u_{n}=  {\sf M}^{-1} T_{n}^* {\sf M}  f.
\label{eq:uv}\end{equation}
Let us show that, for all $g\in L^2 (M)$, 
\begin{equation}
({\sf L} u_{n} ,  g)_{M} =( u_{n}, {\sf L}_{*} g)  .
\label{eq:LMM2a}\end{equation}
 To that end, we use the  Fubini theorem, which requires some estimates on the functions $u_{n}$. They are given in the following assertion.

  \begin{lemma}\label{LTM4}
  If ${\sf M}f\in C_{0}^\infty ({\Bbb R})$, then functions \e{eq:uv} satisfy estimates
    \begin{equation}
| u_{n}  (t)| \leq C_{n}(b) t^b, \q \forall b\in {\Bbb R}.
\label{eq:LMM8}\end{equation}
 \end{lemma}

 \begin{pf}
 Put $\varphi= {\pmb \Gamma} {\sf M} f \in C_{0}^\infty ({\Bbb R})$. Then, by definitions \e{eq:LTM8} and  \e{eq:uv},   $u_{n}= {\sf M}^{-1}  {\pmb \Gamma}^{-1} \hat{\chi}_{n} \varphi$. According to \e{eq:LTM7},   we have  
   \begin{equation}
 u_{n} (t)=   \int_{-\infty}^\inftyÊ G_{n}(t,\eta) \varphi(\eta)   d\eta 
\label{eq:LMM6}\end{equation}
where
  \begin{equation}
G_{n}(t,\eta) = (2\pi )^{-1/2} t^{-1/2}\int_{-\infty}^\inftyÊ Ê t^{i \xi}\Gamma (1/2+i\xi)^{-1}\hat{\chi}_{n}  (\xi -\eta)   d\xi 
\label{eq:LMM6N}\end{equation}
and  $\hat{\chi}_{n}$ is  function \e{eq:LTM7x}.
Here we have used equality \e{eq:LMMv} and the Fubini theorem to interchange integrations over $\xi$ and $\eta $ in the right-hand side of  \e{eq:LMM6}. For an arbitrary $a\in {\Bbb R}$, the integration in \e{eq:LMM6N} can be shifted to the line ${\Bbb R} +ia$  whence
 \begin{equation}
G_{n}(t,\eta) = t^{-a -1/2 } (2\pi )^{-1/2} \int_{-\infty}^\inftyÊ Ê t^{i \xi}\Gamma (1/2-a+i\xi)^{-1}\hat{\chi}_{n}  (\xi +ia-\eta)   d\xi . 
\label{eq:LMM7}\end{equation}
The modulus of the integrand here does not depend on $t$.
In view of the Stirling formula  we have
\[
|\Gamma (1/2-a+i\xi)|= (2\pi)^{1/2} |\xi|^a e^{-\pi |\xi|}\big(1+ O ( |\xi|^{-1})), \q  |\xi|\to \infty.
\]
Therefore it follows from equalities  \e{eq:LTM7x} and  \e{eq:LMM7}  that
  \[
| G_{n}(t,\eta) |\leq C_{n} (a) t^{-a-1/2}\int_{-\infty}^\inftyÊ Ê(1+|\xi|)^{-a} e^{\pi |\xi|} e^{-n^2 (\xi-\eta)^2/4} d\xi. 
\]
The integral here  is bounded by a constant   which depends on $a$ and $n$ but does not depend on $\eta$ in compact intervals. According to  \e{eq:LMM6} this yields estimate \e{eq:LMM8}.
   \end{pf}

 Observe now that under the assumptions of Lemma~\ref{LTM4} the integral 
 \begin{equation}
 (   {\sf L} u_{n}, g)_{M} =\int_{0}^\infty \big(\int_{0}^\infty e^{-\lambda t } u_{n}  (t) dt\big) \ov{g(\lambda)} dM(\lambda)
\label{eq:LMM9}\end{equation}
 converges absolutely because
\[
\int_{0}^\infty e^{-\lambda t }|u_{n}  (t) |dt \leq C (a)(1+\lambda)^{-k/2}, \q \forall k>0,
\]
and, by the Schwarz inequality and condition \e{eq:LTM2},
\[
\int_{0}^\infty (1+\lambda)^{-a} |g(\lambda)| dM(\lambda)\leq \sqrt{\int_{0}^\infty (1+\lambda)^{-k}  dM(\lambda)} \| g\|_{M}<\infty.
\] 
Therefore, by the Fubini theorem, we can interchange    the order of integrations in \e{eq:LMM9} which yields equality \e{eq:LMM2a}.

It follows from   relations
\e{eq:LMM2} -- \e{eq:LMM2a}  that
\[
  (     {\sf L} f, \chi_{n} g)_{M} =(  {\sf M} f,  T_{n}   {\sf M} {\sf L}_{*} g) 
\]
if $ {\sf M} f \in C_{0}^\infty ({\Bbb R}_{+})$ and $g \in L^2 (M)$,   ${\sf L}_{*} g \in L^2 $. Since the operators $T_{n}$ are bounded, this equality    extends to all $f\in L^2  $. If ${\sf L} f \in L^2 (M)$, then using Lemma~\ref{LTM3} we can
pass here to the limit $n\to \infty$ which yields
 relation  \e{eq:LTM3}. 
 
Putting together \e{eq:XL} and \e{eq:LTM3}, we obtain
relation   \e{eq:LTM1} for all $ f \in L^2 $ such that ${\sf L} f \in L^2 (M)$ and all  $g \in L^2 (M)$   such that   ${\sf L}_{*} g \in L^2 $. Hence $ A \subset A_{0}^{**}$.

  Let us summarize the results obtained.
     
       \begin{theorem}\label{QFQL}
  Let $dM(\lambda)$ be a    measure on $[0,\infty)$ satisfying    condition
  \e{eq:Bern1x}. Then the operator $A_{0}:L^2  \to   L^2 (M)$ defined on the domain  ${\cal D}(A_{0})=\cal E$ $($or ${\cal D}(A_{0})= C_{0}^\infty ({\Bbb R}_{+}))$ by the equality $A_{0}f= {\sf L}f$ admits the closure   if and only if $  M(\{0\} )=0$. If $  M(\{0\} )=0$ and  assumption \e{eq:LTM2} is satisfied for some $k>0$, then the  closure $\bar{A}_{0}=:A$ of $A_{0}$   is given   by the same equality $A f={\sf L} f$ on the domain   ${\cal D} (A)$ which  consists of  all $f \in   L^2  $ such that ${\sf L} f\in L^2 (M)$.    
  \end{theorem}
  
    \medskip 
  
   {\bf 3.3.}
  Now we return to Hankel operators. Let us reformulate Theorem~\ref{QFQL} in their terms.
 
 \begin{theorem}\label{QFQ}
  Let $dM(\lambda)$ be a    measure on $[0,\infty)$ satisfying    condition
  \e{eq:Bern1x}. 
Then the form \e{eq:QF} defined on the set $\cal E$ $($or $ C_{0}^\infty ({\Bbb R}_{+}))$ admits the closure  in the space $L^2 ({\Bbb R}_{+})$ if and only if $  M(\{0\} )=0$.  If $  M(\{0\} )=0$ and  assumption \e{eq:LTM2} is satisfied for some $k>0$, then the form \e{eq:QF} is  closed on the set ${\cal D}[h]$ of all  $f \in   L^2  $ such that the integral in the right-hand side of \e{eq:QF}  is finite.    
  \end{theorem}
 
Let us now take for $H$ the self-adjoint nonnegative operator corresponding to the closed form $h[f,f]= {\pmb\la} h, \bar{f}\star f{\pmb\ra}$. It means that 
${\cal D} (H)\subset{\cal D}[h]$ and $(Hf_{1}, f_{2})=h[f_{1}, f_{2}]$ for all $f_{1}\in {\cal D} (H)$ and all
$f_2\in{\cal D} [h]$. Note that ${\cal D}[h]={\cal D}(\sqrt{H})={\cal D} (A)$ and $\| \sqrt{H} f\|= \| A f\|_{M}$ for all $f\in {\cal D} (A)$.

 \begin{remark}\label{QFrem}
 If  $  M(\{0\} )=0$ and  only condition
  \e{eq:Bern1x} is satisfied, then the form \e{eq:QF} is  closed on ${\cal D}[h]= {\cal D} (A_{0}^{**})\subset {\cal D} (A)$. In this case $H$ is still defined as a self-adjoint  operator, but we do not have an explicit description of ${\cal D}(\sqrt{H})$.  
  \end{remark}
  
We can also give a condition on $dM(\lambda)$ guaranteeing that $H>0$.
 
      \begin{proposition}\label{QFP}
  Let   $  M(\{0\} )=0$ and  let assumption \e{eq:LTM2} be satisfied for some $k>0$.  Suppose that every
  set $X\subset {\Bbb R}_{+}$ of full $M$-measure $($that is, $M({\Bbb R}_{+}\setminus X)=0)$ contains infinite number of points $\lambda_{1}, \lambda_{2},\ldots$ such that $\lambda_{n}\to\lambda_{0}>0$ as $n\to\infty$.
 Then  
  $0$ is not an eigenvalue of the operator $H$, that is, $H > 0$.   
  \end{proposition}
        
\begin{pf}
   Indeed, if   $Hf=0$, then  according to \e{eq:QF} we have $({\sf L}f)(\lambda)=0$ for almost all $\lambda\in {\Bbb R}_{+}$ with respect to the measure $M$. It follows that $({\sf L}f)(\lambda_{n})=0$ for some infinite sequence  $\lambda_{n}\to\lambda_{0}>0$. 
      Since the function $({\sf L}f)(\lambda) $ is analytic in the right half-plane, we see that $({\sf L}f)(\lambda) $ for all $\lambda>0$. This implies that  $f=0$ because the kernel of the operator ${\sf L}$ considered in the space $L^2 ({\Bbb R}_{+})$  is trivial.  
\end{pf}

 In view of formula  \e{eq:bbr7} under the assumptions of Theorem~\ref{QFHH}, the measure $dM(\lambda)= \sigma(\lambda) d\lambda$ satisfies the conditions of Proposition~\ref{QFP}. Thus both Theorem~\ref{QFQ} and Proposition~\ref{QFP} are applicable in this case. This yields  all the results stated in   Theorem~\ref{QFHH}.  More generally, it is true for  sums
 \[
 h(t)=\sum_{n=1}^N \kappa_{n} h_{n}(t), \q \kappa_{n}>0,
 \]
 where each function $h_{n}(t)$ satisfies the assumptions of Theorem~\ref{QFHH}.  
  
 Observe that kernels \e{eq:QCs} and \e{eq:QC} may have arbitrary power singularity at the point $t=0$, but it is always required that $h(t)\to 0$ as $t\to \infty$. Without this assumption, there is no reasonable way to define a Hankel operator.
 For example, for $h(t)=1$ representation \e{eq:Bern1} holds with the measure such that $M({\Bbb R}_{+})=0$ and 
 $M(\{ 0\})=1$. Therefore the form \e{eq:QF} does not admit the closure.

   \begin{remark}\label{CC}
 Theorem~\ref{QFQ} can be formulated in a somewhat more general form. Suppose that $h\in C_{0}^\infty ({\Bbb R}_{+})$, and let $\sigma\in{\cal Y}'$ be the corresponding sigma-function. Assume that for all $w\in{\cal Y}$ 
  \begin{equation}
 \la \sigma, w^* w\ra=\int_{0}^\infty | w(\lambda)|^2 dM(\lambda)- \sigma_{0} \| w\|^2_{L^2 ({\Bbb R}_{+})}
\label{eq:Cp}\end{equation}
 where the measure $dM(\lambda)$ satisfies the conditions of  Theorem~\ref{QFQ} and $\sigma_{0} <0$ (if $\sigma_{0} >0$, then \e{eq:Cp} implies that \e{eq:QF}  holds true with the measure  $dM(\lambda)+\sigma_{0}  d\lambda$). 
It follows from \e{eq:Cp}   that the kernel $\ti{h}(t) = h(t)+ \sigma_{0}t^{-1}$ satisfies the assumptions of Theorem~\ref{QFQ} which allows us to define the operator $\wt{H}$ with kernel $\ti{h}(t)$. Since the Carleman operator $\bf C$ is bounded, we can now set $H= \wt{H}-\sigma_{0} \bf C$. 
 \end{remark}

 \medskip 
  
   {\bf 3.4.}
   In the case $\alpha=0$, $r=0$ one can obtain an additional spectral information about the operator $H$. Let us introduce the unitary operator
\[
({\sf D} (\gamma) f)(t)= \gamma^{1/2} f(\gamma t), \q \gamma>0,
\]
of dilations in the space $L^2 ({\Bbb R}_{+})$. The following assertion is intuitively obvious, but it requires a proof because we do not have an explicit description of ${\cal D} (H)$.

 \begin{lemma}\label{UAc}
 If $h(t)=t^{-q}$, then ${\cal D} (H)$ is invariant with respect to ${\sf D}(\gamma)$ and
 \begin{equation}
 {\sf D} (\gamma)^* H {\sf D}(\gamma)=\gamma^{q-1} H.
\label{eq:Ca}\end{equation}
 \end{lemma}
 
 \begin{pf}
 According to \e{eq:bbr7} we now have
$\sigma(\lambda)=   \Gamma (q)^{-1} \lambda^{q-1}$.  
Since ${\sf L} {\sf D} (\gamma) ={\sf D} (\gamma^{-1}){\sf L}  $, it follows  from Theorem~\ref{QFQ}   that 
  ${\cal D}[h] ={\cal D} (\sqrt{H})$ is invariant with respect to the dilations $ {\sf D} (\gamma)$ and
 % \[ h [ {\sf D} (\gamma) f_{1} , {\sf D} (\gamma) f_{2} ]= \gamma^{q-1} h[f_{1},f_{2}].  \]
%  It means that
 \begin{equation}
(\sqrt{H} {\sf D} (\gamma)f_{1},  \sqrt{H} {\sf D}(\gamma)f_{2}) =\gamma^{q-1} (\sqrt{H}  f_{1},  \sqrt{H}  f_{2})
\label{eq:Ca2}\end{equation}
for all $f_{1}, f_{2}\in {\cal D} (\sqrt{H})$.

Let us set $G=\gamma^{(1-q)/2} {\sf D} (\gamma)^* \sqrt{H} {\sf D} (\gamma)$. Then equality \e{eq:Ca2} can be rewritten as
 $(G f_{1}, G f_{2}) =  (\sqrt{H}  f_{1},  \sqrt{H}  f_{2})$.
If $f_{1}\in {\cal D} (H)$, then $  (G f_{1}, G f_{2}) =  (H  f_{1},    f_{2})$ so that $G f_{1} \in {\cal D} (G^*)$  and
$(G^* G f_{1}, f_{2}) =  (H f_{1},      f_{2})$. Since $G=G^*$, it follows that $G^2=H$ which is equivalent to \e{eq:Ca}.
 \end{pf}

Let $X\subset{\Bbb R}_{+}$ be an arbitrary Borel set, and let $E(X)$ be the spectral measure of the operator $H$. Then relation  \e{eq:Ca} can equivalently be rewritten as
 \begin{equation}
 {\sf D}(\gamma)^* E (X) {\sf D} (\gamma)= E (\gamma^{ 1-q}X).
\label{eq:Ca1}\end{equation}

Each of relations \e{eq:Ca} or \e{eq:Ca1} implies that if $\lambda>0$ belongs to the spectrum $\spect (H)$ of the operator $H$, then all points 
$\gamma^{ q -1} \lambda$ also belong to the set $\spect  (H)$. It follows that $\spect  (H)= [0, \infty)$. If $\lambda>0$ is an eigenvalue of $H$, then all points 
$\gamma^{  q-1} \lambda$ are also eigenvalues of $H$. This is impossible so that the operator $H$ does not have eigenvalues. 

Actually, we have a more general statement.\footnote{The author thanks A.~A.~Lodkin and B.~M.~Solomyak for useful consultations on the measure theory.}

 \begin{proposition}\label{UA}
Let $H$ be a   self-adjoint positive operator such that the operators $H$ and $aH$ are unitarily equivalent for all $a>0$. 
Then the spectrum of the operator $H$ coincides with the positive half-line, it has a constant multiplicity  and is absolutely continuous.
 \end{proposition}
 
 \begin{pf}
 According to the spectral theorem we can realize $H$ (see, e.g., the book \cite {BS}) as the operator of multiplication by independent variable $\lambda$ in the space 
 $ L^2({\Bbb R}_{+}; dM(\lambda); {\goth N}(\lambda))$
 where $dM(\lambda)$ is a measure of maximal type with respect to $H$ and 
 $\dim{\goth N}(\lambda)$ equals the multiplicity of the spectrum of the operator $H$  for almost all (with respect to $dM(\lambda)$) $\lambda\in {\Bbb R}_{+}$. Since the operators $H$ and $aH$ are unitarily equivalent, for an arbitrary Borelian set $X\subset {\Bbb R}_{+}$ the conditions $M(X)=0$ and $M(aX)=0$ are equivalent for all $a>0$. This implies that the measure $d M(\lambda)$ is equivalent to the Lebesgue measure on
${\Bbb R}_{+}$. This is proven in Problem~2.12 of Chapter X of  the book \cite{Maka}. To be precise,   the invariance of measures with respect to translations was considered in \cite{Maka}, but the invariance   with respect to dilations  reduces to this case by a change of variables. Thus the operator $H$ is absolutely continuous.

It remains to check that the multiplicity of the spectrum of $H$ is constant. Let $X=X_{k}$ be the Borelian set where this multiplicity   is $k$. Suppose that the Lebesgue measure $|X|>0$. We have to check that $X$ has full measure in ${\Bbb R}_{+}$.
Since the operators $H$ and $aH$ are unitarily equivalent, the sets $X$ and $aX$ coincide up to a set of  the Lebesgue measure zero.   Let $\lambda $ be a density point of $X$, that is
\[
\lim_{\varepsilon\to 0}(2\varepsilon)^{-1} |(\lambda -\varepsilon, \lambda  +\varepsilon)\cap X|=1.
\]
Recall (see, e.g., the book  \cite {Natan}) that almost all points of $X$ possess this property. Evidently, $a\lambda $ is a density point for the set $aX$.
Since the sets $X$ and $aX$ coincide up to a set of  the Lebesgue measure zero, $a\lambda $ is a density point for the set $X$ for all $a >0$. Suppose that $|{\Bbb R}_{+}\setminus X|>0$ and take a density point $\mu\in {\Bbb R}_{+}\setminus X$. Choosing $a= \mu/\lambda$, we see that $\mu=a\lambda$ is also a density point for the set $X$. It follows that
\[
|(\mu -\varepsilon, \mu  +\varepsilon) |=|(\mu -\varepsilon, \mu  +\varepsilon)\cap X|
+|(\mu -\varepsilon, \mu  +\varepsilon)\cap ({\Bbb R}_{+}\setminus X)| =4\varepsilon(1+o(\varepsilon))
\]
while the left-hand side of this relation is $2\varepsilon(1+o(\varepsilon))$.
     \end{pf}

In view of relation \e{eq:Ca}  this result can be directly applied to Hankel operators which yields Theorem~\ref{Car}.
   As shown in the paper \cite{MPT}, the multiplicity of the spectrum of a positive {\it bounded} Hankel operator does not exceed $2$. Most probably, the multiplicity 
 of the spectrum of  the operator $H$ considered in Theorem~\ref{Car} is $1$ because its kernel
 $h(t)=t^{-q}$  has only one singular point $t=\infty$ for $q<1$ and only one singular point $t= 0$ for $q>1$. But this question is out of the scope of the present article.

The results of Theorem~\ref{Car} are of course true for all operators unitarily equivalent to $H$.  For example,  let ${\sf J}$ be the involution in $L^2 ({\Bbb R}_{+})$ defined by the relation $({\sf J}f)(t)=t^{-1} f( t^{-1})$. Then the operator $K={\sf J}^* H {\sf J}$ acts according to the formula
\[
(K f) (t)= \int_{0}^\infty (t+s)^{-q} (ts)^{q-1} f(s) ds.
\]
Therefore all the conclusions of Theorem~\ref{Car} are satisfied  for such operators $K$ (if $q\neq 1$).

 \medskip 
  
   {\bf 3.5.}
    In the general case some spectral information is also available. Below we consider the quadratic form of $H$ on the characteristic functions $\mathbbm{1}_{(a,b)}(t)$ of intervals $(a, b)\subset{\Bbb R}_{+}$. Obviously, $\mathbbm{1}_{(a,b)}\in{\cal E}$ if $ 0<a<b<\infty$, $\| \mathbbm{1}_{(a,b)}\|^2= b-a$ and
     \begin{equation}
({\sf L} \mathbbm{1}_{(a,b)})(\lambda)=  (e^{-a\lambda}-e^{-b\lambda})\lambda^{-1}
\label{eq:ch}\end{equation}
    
      \begin{proposition}\label{GSI1}
If  $  M(\{0\} )=0$ and   condition
  \e{eq:Bern1x} is satisfied, then  the point zero belongs to   the spectrum of the corresponding Hankel operator $H$. \end{proposition}
 
 \begin{pf}
 Let $f_{n}= \mathbbm{1}_{(n,n+1)}$.
According to \e{eq:ch} the functions $({\sf L} f_{n})(\lambda) $ are uniformly bounded by $e^{-\lambda}$ and tend to zero  as $n\to \infty$ for all $\lambda>0$. Hence, by the dominated convergence theorem, it follows from \e{eq:QF}  that $\|\sqrt{H}f_{n}\|\to 0$ as $n\to\infty$. Thus $0\in \spec (H)$.
   \end{pf}
   
   Of course this  result is consistent with the general fact (see, e.g., the book \cite{Pe}) that  $0\in \spec (H)$ for all {\it bounded} Hankel operators.
   
      \begin{proposition}\label{GSI2}
Let $  M(\{0\} )=0$, and   let condition
  \e{eq:Bern1x} be satisfied. If at least one of the  conditions   \e{eq:Wid} is violated,  then  the corresponding Hankel operator $H$ is unbounded. 
  \end{proposition}
 
 \begin{pf}
 If the first condition  \e{eq:Wid} is violated,  then there exists a sequence $\varepsilon_{n}\to 0$ such that $\varepsilon_{n}^{-1} M(0,\varepsilon_{n}) \to\infty$ as $n\to\infty$. Put $f_{n}= \mathbbm{1}_{(1, \varepsilon_{n}^{-1} )}$. It follows from \e{eq:ch} that $|({\sf L} f_{n})(\lambda)| \geq (e\varepsilon_{n})^{-1}$   for $\lambda\in (0, \varepsilon_{n})$. Hence according to \e{eq:QF} we have 
$
 \|\sqrt{H}f_{n}\|^2\geq  (e \varepsilon_{n})^{-2} M(0,\varepsilon_{n}) 
$
 so that $\| f_{n}\|^{-2} \|\sqrt{H}f_{n}\|^2\geq e^{-2} \varepsilon_{n}^{-1} M(0,\varepsilon_{n})  \to \infty$ as $n\to\infty$.
 
  Similarly, if the second condition  \e{eq:Wid} is violated,  then there exists a sequence $l_{n}\to \infty$ such that $l_{n}^{-1} M(0, l_{n}) \to\infty$ as $n\to\infty$.  Put $f_{n}= \mathbbm{1}_{( l_{n}^{-2} , l_{n}^{-1} )}$.
    It follows from \e{eq:ch} that $|({\sf L} f_{n})(\lambda)| \geq (e l_{n})^{-1}$   for $\lambda\in (0, l_{n})$. Hence according to \e{eq:QF} we have 
 $
 \|\sqrt{H}f_{n}\|^2\geq ( e  l_{n})^{-2} M(0, l_{n}) 
$
 so that again $\| f_{n}\|^{-2} \|\sqrt{H}f_{n}\|^2\geq e^{-2} l_{n}^{-1} M(0, l_{n}) \to \infty$ as $n\to\infty$.
  \end{pf}
  
  In view of formula \e{eq:bbr7} Propositions~\ref{GSI1} and \ref{GSI2} directly apply to
  Hankel operators $H$ with kernels \e{eq:QCs} and \e{eq:QC}.  
  
   Proposition~\ref{GSI2} is essentially equivalent to the result of H.~Widom mentioned in Section~1, but our proof relies on the construction of trial functions and is quite different from that in \cite{Widom}.
  
      %  Non-compact perturbations.

%************************************************************
  %************************************************************
\section{Perturbation theory}  
 %************************************************************
 
  Here we study    perturbations of singular quasi-Carleman operators   $H_{0}$ introduced in Section~3 by  bounded and, in particular, compact   self-adjoint Hankel operators $V$ with kernels \e{eq:E1v}. 
  
  %  Our results imply that, under very general assumptions,  the operators $H=H_{0}+V$ and $V$ have the same numbers of negative eigenvalues.

  \medskip
    
  {\bf 4.1.}
   As far as the unperturbed operator $H_{0}$ is concerned, we accept the following

   \begin{assumption}\label{Hfree}
The sigma-function $\sigma_{0}(\lambda) $ of the operator $H_{0}$ is nonnegative, $\supp \sigma_{0}\subset [0,\infty)$,      $\sigma_{0}\in L^\infty_{\rm loc} ({\Bbb R}_{+})$   and 
  \begin{equation}
 \sigma_{0}(\lambda)= O(\lambda^{- l_{+}}) \;{\rm as} \;  \lambda\to 0,   \q 
  \sigma_{0}(\lambda)= O(\lambda^{l_{-}})  \; {\rm as} \;  \lambda\to \infty 
  \label{eq:Expsi}\end{equation}
where $l_{+} <1$ and $l_-$ may be arbitrary large.
  \end{assumption}

  This assumption can of course be equivalently reformulated in terms of the sign-function $s_{0}(x)=\sigma_{0} (e^{-x})$: 
  $s_{0}(x) \geq 0$, $\supp s_{0}\subset {\Bbb R}$,      $s_{0}\in L^\infty_{\rm loc} ({\Bbb R} )$   and 
 \[
 s_{0}(x)= O(e^{l_{+}x})  \;{\rm as} \;  x\to +\infty \q{\rm and} \q s_{0}(x)= O(e^{l_{-}|x|})  \;{\rm as} \;  x\to -\infty .
 \]
 Of course, Assumption~\ref{Hfree} does not guarantee that $s_{0} \in{\cal S}'$.

  Since the measure $dM_{0}(\lambda)=\sigma_{0}(\lambda) d\lambda$ satisfies condition \e{eq:LTM2},   according to Theorem~\ref{QFQ}   the form 
  \begin{equation}
h_{0}[f,f]  =  \int_{0}^\infty   | ({\sf L}f)(\lambda) |^2  \sigma_{0}(\lambda) d\lambda
\label{eq:C1Mm}\end{equation}
is closed on the set ${\cal D}[h_{0}]$ of all functions $f\in L^2 ({\Bbb R}_{+})$ such that integral \e{eq:C1Mm} is finite.     Setting $\lambda=e^{-x}$, we  can equivalently   rewrite definition \e{eq:C1Mm}
   as
       \begin{equation}
h_{0}[f,f]  =  \int_{-\infty}^\infty   |u(x) |^2  s_{0}(x) dx
\label{eq:CB}\end{equation}
where  $u(x)$ and $f(t)$ are linked by formula \e{eq:sign2}.
 We denote by 
$H_{0}$ the  self-adjoint operator corresponding to the form \e{eq:C1Mm} (or \e{eq:CB}). Of course, $H_{0}\geq 0$. Moreover, by Propositions~\ref{QFP} and \ref{GSI1}, the point $0\in\spec (H_{0})$, but it is not the eigenvalue of $H_{0}$.

 It follows from  formula \e{eq:bbr7}   that  Assumption~\ref{Hfree} is satisfied for kernels $h_{0}(t)$  given by equality  \e{eq:E1r} where either $\alpha=0$ and $q>0$ or $\alpha >0$ and $q\geq 1$. Now $l_{+}= 1-q$ for $\alpha=0$ and 
 $l_{+} $ is arbitrary for $\alpha>0$; $l_{-}=  q-1$ for $r=0$ and 
 $l_{-} $ is arbitrary for $r>0$. 
  In the   case $\alpha=0$,  $q>0$  the operators $H_{0}$ are unbounded unless $q=1$ when $H_{0}=\bf C$ is the Carleman operator. If $\alpha >0$, but $r=0$, then the operators $H_{0}$ are unbounded   for $q> 1$, but $H_{0}$ is  bounded   for $q= 1$. If $\alpha >0$ and $r>0$, then the operators $H_{0}$ are compact for all values of $q$.

   Another interesting example are Hankel operators  $H_{0}$ with kernels
   \begin{equation}
h_{0}(t)= P (\ln t)t^{-1}
\label{eq:LOG}\end{equation}
where  
  $
P (x)=\sum_{k=0}^K p_{k} x^k$,  $ p_{K}>0 $,
 is an arbitrary real polynomial of {\it even} degree $K$. Such operators were studied in \cite{Yf1} where, in particular, it was shown that $H_{0}$ are semibounded and the essential spectrum $\spec_{\rm ess} (H_{0})=[0,\infty)$ unless $K=0$.
 The sign-function of kernel \e{eq:LOG} is given by the polynomial
$
   s_{0}(x)=  \sum_{ k=0}^K q_{k} x^k 
$
 where the coefficients
$q_{k}$ admit an explicit expression in terms of the coefficients $p_{k}, p_{k+1}, \ldots, p_{K}$.
 Of course  $s_{0} \in {\cal S}'$ and Assumption~\ref{Hfree} is satisfied for kernels \e{eq:LOG} provided $s_{0}(x)\geq  0$. 
 Note that the inequality $ s_{0}(x)\geq 0$ implies that $P(x)\geq 0$ but not vice versa. 
 
On the contrary, sigma-functions of finite rank Hankel operators are singular distributions (see Section~5) so that
  Assumption~\ref{Hfree} is violated for  such operators.
        
\medskip
    
  {\bf 4.2.}
  Let us start with a general statement on perturbations of operators $H_{0}$ satisfying Assumption~\ref{Hfree} by bounded Hankel operators $V$. Recall that, for all $f\in L^2 ({\Bbb R}_{+})$, we have
  \[
  (Vf,f) ={\pmb \la} \sigma_{v}, w^* w {\pmb \ra} =: \sigma_{v}[w,w]
  \]
  where $w={\sf L}f\in {\Bbb H}_{r}^2$ and $\sigma_{v}\in ({\Bbb H}_{r}^1)'$.
   We put $H =H_{0} + V$.  This operator is defined via its quadratic form
  \begin{align}
  h[f,f]=& h_{0}[f,f]+ (Vf,f) 
   \nonumber\\
 =&  \int_{0}^\infty \sigma_{0} (\lambda) |w(\lambda)|^2 d\lambda + \sigma_{v}[w,w]= : \sigma [w,w], \q w={\sf L}f.  
   \label{eq:DHW}\end{align}
    This form    is closed on the set ${\cal D}[h]$ of
     all $f\in L^2 ({\Bbb R}_{+})$ (or equivalently of all $w\in {\Bbb H}_{r}^2$) such that  the integral in the right-hand side converges. Of course ${\cal D}[h] ={\cal D}[h_{0}]$. The following auxiliary result (cf. Theorem~\ref{HBx}) is a direct consequence of equality \e{eq:DHW}.
  
  \begin{lemma}\label{GGL}
 Let   Assumption~\ref{Hfree} be  satisfied, and let    a Hankel operator $V$ be bounded. Then $N_{-} (H)$ equals the maximal dimension of linear sets ${\cal K}\subset{\Bbb H}_{r}^2$ such that
 $ \sigma [w,w]<0$
    for all $w\in {\cal K}$, $w\neq 0$.
    \end{lemma}
    
    Note  that if $w\in {\Bbb H}_{r}^2$ and  $ \sigma [w,w]<0$, then  $ \sigma _{0}[w,w]< \infty$ and hence $f={\sf L}^{-1}w\in {\cal D}[h]$.

If $s_{0}\not\in {\cal S}'$, 
we cannot use Theorem~\ref{HByz} and  study the form \e{eq:DHW} on  the set $C_{0}^\infty ({\Bbb R}_{+})$ of test functions $w(\lambda)$
(this reduction was essentially used by the proof of Theorem~\ref{HKL} in \cite{Yqd}).   Now we are obliged to work with   analytic functions $w\in {\Bbb H}_{r}^2$ (or the corresponding  test functions $u(x)= e^{-x/2}w(e^{-x})$) which  is technically more involved. In the case 
    $s_{0}\in{\cal S}'$ (for example, for Hankel operators with kernels \e{eq:LOG}) the proofs below can be considerably simplified. On the other hand, the advantage of the method suggested here is that  
    it  directly yields, by  formula   \e{eq:sign2}, trial functions $f(t)$ for which $h[f,f]<0$. 
    
    % Apparently, under the assumptions of Theorem~\ref{FDH} below it is hardly possible to guess such functions $f(t)$ in the original representation $L^2 ({\Bbb R}_{+})$.
    
In many cases it suffices to  consider   gaussian trial functions
    \begin{equation}
  u_\varepsilon(x; a)=  \varepsilon^{-1/2}e^{-(x-a)^2/\varepsilon^2}, \q \varepsilon>0,\q a\in{\Bbb R}.
  \label{eq:GF}\end{equation}
  Obviously, $\|u_\varepsilon(a)\|^2= \sqrt{\pi/2} $. Since
   \begin{equation}
  (\Phi  u_\varepsilon( a) )(\xi)= 2^{-1/2}\varepsilon^{1/2}e^{-i\xi a}e^{-\varepsilon^2 \xi^2/4},
  \label{eq:GF1}\end{equation}
  functions $f_\varepsilon(t; a)$ defined by \e{eq:sign2}  
  belong to $L^2 ({\Bbb R}_{+})$ or, equivalently, the corresponding functions  $w_\varepsilon(\lambda; a)$ 
   belong to ${\Bbb H}^2_{r}$. Note however that $\|f_\varepsilon(a)\|\to \infty$ as $\varepsilon\to 0$.
  
  The following assertion is almost obvious.

   \begin{lemma}\label{GF}
If the parameters $a_{1},\ldots, a_{N}$ are pairwise different, then the functions $ u_\varepsilon( x; a_{1}),\ldots, u_\varepsilon( x; a_{N})$ are linearly independent.
  \end{lemma}
  
   \begin{pf}
   It suffices to check that the functions $ (\Phi u_\varepsilon( a_{1}))(\xi),\ldots, (\Phi u_\varepsilon( a_{N}))(\xi)$ are linearly independent or according to \e{eq:GF1} that the functions $e^{-i\xi a_{1}},\ldots, e^{-i\xi a_{N}}$ are linearly independent. If $\sum_{j=1}^N c_{j}  e^{-i\xi a_j} =0$, then differentiating this identity and putting $\xi=0$, we obtain the system 
 $\sum_{j=1}^N c_{j}    a_j^n =0$  where $n=0,1,\ldots, N-1$  for $c_{1},\ldots, c_{N}$. The determinant of this system  is the Vandermonde determinant. Since it is not zero, we see that $c_{1}=\cdots =c_{N}=0$. 
   \end{pf}

% \medskip  {\bf 4.3.}

 For Hankel operators $V$ with regular sign-functions, we use the following assertion.
  
   \begin{theorem}\label{GG}
 Let  Assumption~\ref{Hfree} be satisfied, and let $H_{0}$ be the corresponding Hankel operator. Suppose that a Hankel operator $V$ is bounded and that its sign-function $s_{v}\in L^1_{\rm loc} ({\Bbb R})$.
 Put $H =H_{0} + V$ and $s=s_{0}+ s_{v}$. 
  Then:
 
   $1^0$ The operator  $H \geq 0$    if $s(x)\geq 0$ for almost all $x\in {\Bbb R}$.
   
      $2^0$  The operator  $H  $  has infinite negative spectrum if $s(x)\leq -s_{0} < 0$ for almost all $x $ in some interval $\Delta\subset  {\Bbb R}$ and $s(x)$ is exponentially bounded away from $\Delta$. 
          \end{theorem}

  \begin{pf}
 The assertion  $1^0$ is obvious because (cf. formula  \e{eq:CB} for $ h_{0}[f,f]$)  for all $f\in {\cal D}[h] $ we have
 \begin{equation}
  h[f,f]= \int_{-\infty}^\infty s(x) |u(x)|^2 dx 
 \label{eq:gB}\end{equation}
where $u(x)$ and $f(t)$ are related by  \e{eq:sign2}.

Let us prove $2^0$.  For an arbitrary $N$, let us choose points $a_{1}, \ldots, a_{N}\in\Delta = :(a_{0}, a_{N+1})$ in such a way that $ a_{j+1}-a_{j}=a_{j }-a_{j-1}=:\d$
  for all $j=1,2, \ldots, N $ and define functions $u_{\varepsilon} (x;a_{j})$ by formula \e{eq:GF}
where $\varepsilon$ is a sufficiently  small number. By Lemma~\ref{GF} the functions $u_{\varepsilon} (x;a_{1}), \ldots, u_{\varepsilon} (x;a_{N})$ are linearly independent.   By our condition on $s(x)$, we have
 \begin{multline*}
\int_{\Delta} s(x) u_\varepsilon^2(x; a_{j})  dx \leq -s_{0}\varepsilon^{-1}\int_{ a_{0}}^{ a_{N+1}} e^{-2(x-a_{j})^2/\varepsilon^2} dx
\\
=-  s_{0} \int_{(a_{0}-a_{j})/\varepsilon}^{ (a_{N+1}-a_{j})/\varepsilon} e^{-2y^2 } dy =-  s_{0} (\sqrt{\pi/2}+O(\varepsilon))\leq -  s_{0} 
 \end{multline*}
because $(a_{0}-a_{j})/\varepsilon\to -\infty$ and $(a_{N+1}-a_{j})/\varepsilon\to \infty$ as $\varepsilon\to 0$.
Moreover,     for all $\d>0$ and a  sufficiently large $l>0$, 
 \begin{equation}
 \int_{{\Bbb R}\setminus\Delta} s(x) u_\varepsilon^2(x; a_{j}) dx
 \leq  C \varepsilon^{-1} \int_\d^\infty  e^{l x} e^{-2x^2/\varepsilon^2} dx=
C  \int_{ \d/\varepsilon}^\infty e^{\varepsilon l y } e^{-2y^2 }  dy =O(\varepsilon) .
 \label{eq:gg2}\end{equation}

  For $j\neq k$, we have
 \[
\big| \int_{\Delta } s(x) u_\varepsilon(x;a_{j}) u_\varepsilon(x; a_{k})dx \big| \leq C \varepsilon^{-1} 
e^{-\d^2/2\varepsilon^2} \int_{\Delta } |s(x) | dx =O(\varepsilon),
\]
and according to \e{eq:gg2}
 the corresponding integral over ${\Bbb R}\setminus\Delta $ is also $O(\varepsilon)$.

 Putting together the estimates obtained, we see that
 \[
 s[\sum_{j=1}^N \nu_{j} u_\varepsilon(a_{j}) , \sum_{j=1}^N \nu_{j} u_\varepsilon(a_{j})]
 \leq - s_{0}\sum_{j=1}^N |\nu_{j}|^2 (1+O(\varepsilon)).
 \]
 Therefore the operator $H$ has at least $N$ negative eigenvalues.  
   \end{pf}

   \begin{remark}\label{GGr}
Under the assumptions of Theorem~\ref{GG} the operator $V$ is not supposed to be compact, but if it is compact then,
 in the assertion  $2^0$,  the negative spectrum of the  operator  $H  $ consists of infinite number of   eigenvalues.
    \end{remark}

  %$\sigma_{\rm ess} (H)=\sigma_{\rm ess} (H_{0})\subset [0,\infty)$ and 
  
 \medskip
    
  {\bf 4.3.}
   Let us now consider perturbations of singular operators $H_{0}$ defined in Theorem~\ref{QFQ}  by Hankel operators $V$ with kernels \e{eq:E1v}.
According to formula \e{eq:bbr7} 
the   corresponding sigma-function $\sigma_{v} (\lambda)$   is given by the equality  
\begin{equation}
\sigma_{v} (\lambda)= v_{0}
   \frac{e^{\beta \rho}}{\Gamma(-k)}   (\lambda - \beta)_{+}^{-k-1}e^{-\rho \lambda} .
\label{eq:bvr7}\end{equation}
We impose  conditions on the pararameters $\beta$, $\rho$ and $k$ such  that the operators $V$ are  bounded.

It turns out that the cases $k<0$ and $k\geq 0$ are qualitatively different. In the first case the sign-function of $V$ is regular, so that the negative spectrum of $H=H_{0}+V$ is governed by Theorem~\ref{GG}. In the second case the negative spectrum of $H=H_{0}+V$ is determined solely by the singularity of the sigma-function $\sigma_{v} (\lambda)$ at the point $\lambda=\beta$.

   Consider first perturbations \e{eq:E1v} for $k<0$.  If $v_{0}\geq 0$, then, by Theorem~\ref{HKL}, the operator $V\geq 0$ and hence $H=H_{0}+ V \geq 0$. So we suppose that $v_{0} <0$.
 If  $k\in (-1,0)$ and $\beta >0$, then $  \sigma_{v} (\lambda)$ is continuous for $\lambda> \beta$,  $ \sigma_{v}  (\lambda)\to -\infty$ and hence $\sigma (\lambda)= \sigma_0 (\lambda)+ \sigma_v (\lambda) \to-\infty$ as $\lambda\to  \beta +0$. 
 Since $ \sigma_{v}\in L^1_{\rm loc} ({\Bbb R}_{+})$,  it follows   from Theorem~\ref{GG} that the   operator $H $ has infinite negative spectrum.

If  $k = -1$ or $k <-1$ but $\rho>0$, then  $\sigma_{v} (\lambda)$
is a bounded negative function. By Theorem~\ref{GG}, the operator $H \geq 0$ if and only if 
$\sigma  (\lambda)\geq 0$ (for all $ \lambda\geq \beta$).
 In the opposite case it has   infinite negative spectrum.

 % for all $\lambda>0$ (we have set $\lambda=e^{-x}$ in \e{eq:bbr}). 

Let us summarize the results obtained.

 \begin{theorem}\label{HKC}
 Let $H = H_{0} + V$ where the sigma-function  $\sigma_{0}(\lambda)$ of the operator $H_{0}$ satisfies Assumption~\ref{Hfree}, and let $V$ be the Hankel operator    with kernel \e{eq:E1v} where $k<0$. Then:
 
   $1^0$ If $k >-1 $ and $\beta>0$, then the operator  $H > 0$ for $v_{0} \geq 0$  and $H $     has infinite negative spectrum for all $v_{0} <0$.
 
  $2^0$  Let $k\leq -1 $. If $k < -1 $   suppose that $\rho>0$. Put
  \[
  \nu = \Gamma(-k) e^{- \beta \rho}\ess\inf_{\lambda\geq \beta} \Big(   (\lambda - \beta)^{k +1}e^{\rho \lambda} \sigma_{0} (\lambda)\Big).
  \]
  Then  the operator  $H \geq 0$     for   $v_{0} \geq-\nu $, and it   has   infinite negative spectrum for  $v_{0} < -\nu $.
    \end{theorem}

For the Carleman operator $H_{0}=\bf C$, we have $\sigma_{0} (\lambda)=1$ and hence $ \nu  = 1$ for all values of $\beta$ and $\rho$ if $k=-1$ and
 \[
  \nu  = \Gamma(-k) e^{-k-1} \big(\frac{\rho}{-k-1}\big)^{-k-1}, \q k < -1.
 \]
In particular,  the critical coupling constant $\nu $   does not depend on $\beta$ in this case.

     \begin{remark}\label{fDH}
In part $1^0$ of Theorem~\ref{HKC} the operator $V$ is compact so that for $v_{0} <0$ the operator $H$ has infinite number of negative eigenvalues accumulating to the point zero. The same is true in part $2^0$ if $v_{0} < -\nu$ and either $\beta>0$ or $\beta=0$, $k<-1$ (in these cases $V$ is compact). If $V$ is not compact, then the operator $H=H_{0} +V$ may have negative continuous spectrum. For example, if $h_{0}(t)= t^{-1}$ and $v(t)= v_{0}t^{-1}$, then the spectrum of $H=(1+v_{0}){\bf C}$ coincides with the interval $[0, (1+v_{0})\pi]$ for $v_{0}>-1$ and with the interval $[  (1+v_{0})\pi,0]$ for $v_{0}< -1$.
 \end{remark} 
 
   The case $k>0$ when the operator $V$ is not sign-definite is essentially more difficult.    Our goal is to prove the following result.
     
      \begin{theorem}\label{FDH}
 Let $H = H_{0} + V$ where the sigma-function  $\sigma_{0}(\lambda)$ of the operator $H_{0}$ satisfies Assumption~\ref{Hfree} and $V$ is the Hankel operator    with kernel \e{eq:E1v}.
  Suppose that $\beta >0$, $\rho\geq 0$ and $k > 0 $
for some $k\not\in {\Bbb Z}_+$. Then:

$1^0$ If $v_{0}>0$ and $[k]$ is odd, then   $N_- (H) = ([k]+ 1)/2$. 

$2^0$ If $v_{0} < 0$ and $[k]$ is even, then $N_- (H )= [k]/2+1$.

 $3^0$ If $v_{0}>0$ and $[k]$ is even or $v_{0}< 0$ and $[k]$ is odd, then $N_- (H )=  \infty$.
 \end{theorem} 
 
Putting together Theorems~\ref{HKL} and \ref{FDH}, we obtain relation \e{eq:FDb}. This implies
    Theorem~\ref{FDHC} for $k\not\in {\Bbb Z}_{+}$. The case $k\in {\Bbb Z}_{+}$ will be considered in the next section.

 % A difficulty of the proof of Theorem~\ref{FDH}  is that $\sigma_{0}\not\in {\cal S}'$ and hence $\sigma\not\in {\cal S}'$. Therefore Theorem~\ref{HBy} (??) is not now directly applicable although its statement remains true. -- to remove; already explained. 

    \medskip
    
  {\bf 4.4.} 
  Let us prove parts $1^0$ and $2^0$ of Theorem~\ref{FDH}. 
   Set $n=[k]$ and    $\ell= [n/2] +1$.
    Since $H_{0}\geq 0$, we have $N_- (H ) \leq N_{-} (V)$, 
  which in view of Theorem~\ref{HKL} yields the upper estimate $N_- (H )\leq \ell $.   According to  Lemma~\ref{GGL}, to prove the lower  estimate $N_- (H )\geq \ell $,  we have to construct a linear subspace 
  ${\cal K}\subset{\Bbb H}^2_{r} $ of   dimension  $\ell$    such that   $\sigma[w,w]<0$ for all $w\in  {\cal K}$, $w\neq 0$.  
  Note that our assumptions $v_{0} \Gamma(-k)^{-1} >0$.
  
  We   need an   elementary assertion about distributions $\mu^{-k-1}_+ $.

    \begin{lemma}\label{SiD}
     Suppose that   $k \in {\Bbb R}_{+} \setminus {\Bbb Z}_{+} $.
Let a bounded $C^\infty$ function $\varphi(\mu)$ of $\mu\in \Bbb R$ satisfy the conditions
 \begin{equation}
\varphi(0) = 1 \q  \mbox{and}\q \varphi' (0) = \cdots = \varphi^{(n)} (0) = 0  \q  \mbox{if} \q n  \geq 1,
\label{eq:L6}\end{equation}
and let $Q(\mu)$ be a polynomial of $\deg Q \leq n$. Then
 \begin{equation}
   \int_{-\infty}^\infty \mu^{-k-1}_+ Q(\mu) \varphi^2 (\mu) d \mu
   =   \int_{0}^\infty \mu^{-k-1} Q(\mu) ( \varphi^2 (\mu) -1) d\mu. 
   \label{eq:did}\end{equation}
 \end{lemma}
 
 \begin{pf}
 According to \e{eq:L6} for the function $\psi(\mu) =Q(\mu ) \varphi^2 (\mu)$, we have
$\psi^{(p)} (0) =  Q^{(p)} (0)$ for all $ p=0,\ldots, n$,
whence
\[
 \sum_{p=0}^n  \frac{1}{p!   } \psi^{(p)}  (0) \mu^p = 
 \sum_{p=0}^n  \frac{1}{p!   } Q^{(p)} (0) \mu^p =  Q(\mu)
 \]
 if $n\geq \deg Q$. Therefore relation \e{eq:did} is a direct consequence of definition  \e{eq:di}.
   \end{pf}

  Put  
      \begin{equation}
w_{\varepsilon} (\lambda) =   P (\lambda- \beta) R (\lambda-\beta ) \exp\big(-\varepsilon^{-2m} \ln^{2m}(\lambda/ \beta)\big)
\label{eq:L5PT}\end{equation}
 where $P  (\mu )= \sum_{j=0}^{\ell -1} p_{j}  \mu ^j$   is an arbitrary polynomial of $\deg P \leq \ell -1$
 and $R (\mu) =\sum_{j=0}^n r_{j}  \mu ^j $ is a special
  polynomial of   $\deg R \leq n$,  $\varepsilon$ is a small parameter and $m$ is a sufficiently large number.  It is easy to see that functions \e{eq:L5PT} belong to the Hardy space ${\Bbb H}^2_{r} $ for all $m=1,2, \ldots$.

 First we construct the polynomial 
 $R (\mu) $.
We require that the function
  \begin{equation}
\theta(\mu)= R(\mu)  e^{- \rho\mu/2}
\label{eq:X5PT}\end{equation}
  satisfy the equations 
   \begin{equation}
\theta(0) = 1 \q  \mbox{and}\q \theta' (0) = \cdots = \theta^{(n)} (0) = 0  \q  \mbox{if} \q n  \geq 1.
\label{eq:L6z}\end{equation}
Solving these equations for $r_{0},\ldots, r_{n}$, we find successively all the coefficients $r_{0}=1, r_{1}=\rho/2, r_2,\ldots, r_{n}$.  Note that in the case $\rho=0$, we have $R(\mu)=1$ so that
 this construction is not necessary.
 
 Let us estimate the sigma-form of the operator $V$.
 
   \begin{lemma}\label{FDnt}
 Let $\sigma_{v}(\lambda)$ be function \e{eq:bvr7}, and let the functions $w_{\varepsilon}$ be defined by equality \e{eq:L5PT} where $2m>n$. Suppose that function \e{eq:X5PT} satisfies conditions \e{eq:L6z}. Then there exist $\varepsilon_{0}>0$ and $c>0$ such that 
   \begin{equation}
v_{0}^{-1}\Gamma(-k)  \sigma_{v} [w_{\varepsilon},w_{\varepsilon}] \leq - c\|P\|^2 ,\q  \|P\|^2 := \sum_{j=0}^{\ell -1}| p_{j}|^2 ,\q \forall \varepsilon\in (0,\varepsilon_{0}).
 \label{eq:L5PS}\end{equation}
 \end{lemma} 
 
   \begin{pf}  
   Put
   \[
\varphi_{\varepsilon} (\mu)=\theta(\mu)\exp\big(-\varepsilon^{-2m} \ln^{2m}(1 +\mu/\beta)\big)
\]
where $\theta(\mu)$ is function \e{eq:X5PT}. 
   According to \e{eq:bvr7} and \e{eq:L5PT} we have the expression 
   \begin{equation}
\sigma_{v}[ w_{\varepsilon},w_{\varepsilon}]=     \frac{v_{0}}{\Gamma(-k)}   \int_{-\infty}^\infty \mu_{+}^{-k-1}|P(\mu)|^2 \varphi_{\varepsilon}^2 (\mu)   d \mu .
\label{eq:L4}\end{equation}
Since $2m>n$,
   equations \e{eq:L6z} for $\theta (\mu)$ imply equations \e{eq:L6} for    
$\varphi_{\varepsilon} (\mu)$. 
Therefore Lemma~\ref{SiD} applied to $Q(\mu) = |P(\mu)|^2$ yields the representation
  \begin{equation}
v_{0}^{-1}\Gamma(-k)  \sigma_{v} [w_{\varepsilon},w_{\varepsilon}] 
=      \int_{0}^\infty |P (\mu)|^2
  \big(\theta^2(\mu)  e^{- 2 \varepsilon^{-2m} \ln^{2m}(1+\mu/ \beta)}- 1 \big) \mu^{-k-1}  d\mu  .
  \label{eq:L7PT}\end{equation}
  Note  that
   \begin{equation}
    \int_{0}^\infty |P (\mu)|^2
  \big(   e^{- 2 \varepsilon^{-2m} \ln^{2m}(1+\mu/ \beta)}- 1 \big) \mu^{-k-1}  d\mu
  \leq      -c \| P\|^2
  \label{eq:hh2}\end{equation}
  for all $\varepsilon\leq 1$. Indeed, the left-hand here is maximal for $\varepsilon=1$ and hence it suffices to use that
  \[
  \min_{\| P \| =1} \int_{0}^\infty |P (\mu)|^2 \big( 1-  e^{- 2   \ln^{2m}(1+\mu/ \beta)} \big) \mu^{-k-1}  d\mu \geq c>0.
  \]
  
If $\rho=0$, then  $\theta(\mu) =1$, and hence estimates \e{eq:hh2} and \e{eq:L5PS} coincide.  If $\rho>0$, we have to get rid of $\theta(\mu) $ in the right-hand side of \e{eq:L7PT}.
   It follows from conditions \e{eq:L6z}    that
$
\big| \theta^2(\mu)-1 | \mu^{-k-1 }\leq C \mu^{n -k }.
$
Using also an obvious  estimate
 \begin{equation}
|P (\mu)|^2 \leq C \| P\|^2 (1+ \mu^{2\ell -2}), 
 \label{eq:hp2}\end{equation}
we see that  
  \begin{multline}
   \int_{0}^\infty |P (\mu)|^2
   |\theta^2(\mu)-1 |  e^{- 2 \varepsilon^{-2m} \ln^{2m}(1+\mu/ \beta)}  \mu^{-k-1}  d\mu 
  \\
  \leq    C \| P\|^2  
     \int_{0}^\infty  
    e^{- 2 \varepsilon^{-2m} \ln^{2m}(1+\mu/ \beta)} ( \mu^{n -k}+  \mu^{n -k+ 2\ell -2})  d\mu .
  \label{eq:hh3}\end{multline}
  Since $n -k>-1$ and $\ell\geq 1$, the integral in the right-hand side of \e{eq:hh3} tends to zero as $\varepsilon\to 0$, and hence  the right-hand side of  \e{eq:hh3} is bounded by
 $c \| P\|^2  /2 $ for sufficiently small $\varepsilon>0$.
  Putting together this result with \e{eq:hh2}, we obtain estimate \e{eq:L5PS}.
 \end{pf}

  Let us now consider the sigma-form $\sigma_{0}[w_{\varepsilon}, w_{\varepsilon}]$ on functions $w_{\varepsilon}$ defined by formula \e{eq:L5PT}. 
  %Note an obvious estimate \[ |P(\lambda-\beta) R (\lambda-\beta)|\leq C \| P\| (1+\lambda^{\ell +n-1}).  \]
Using again estimate \e{eq:hp2}, conditions \e{eq:Expsi} where $l_{+}  <1$ and making the change of variables $\lambda= \beta e^{\varepsilon x}$, we see 
      that
     \[
\sigma_0[w_{\varepsilon},w_{\varepsilon}] 
\leq  C \|P\| ^2\int_{0}^\infty  (\lambda^{-l_{+}}  
+ \lambda^{ 2(\ell+n-1) +l_{-}} )  \exp\big( -2 \varepsilon^{-2m} \ln^{2m}(\lambda/ \beta)\big) d\lambda
\leq  C_{1} \varepsilon \|P\| ^2.
\]

 Putting together this  estimates  with \e{eq:L5PS}, we see that
$\sigma [w_{\varepsilon},w_{\varepsilon}] <0$
if $w_{\varepsilon}$  is defined by formula \e{eq:L5PT} where $\varepsilon$ is sufficiently small and  $P (\mu) $ is an arbitrary nontrivial polynomial of $\deg P \leq \ell-1$. 
Since the dimension of such polynomials equals $\ell$, this yields us the 
linear subspace ${\cal K} \subset {\Bbb H}^2_{r} $ of   dimension $\ell$ where the form $\sigma$ is  negative.  This shows that $N_{-} (H)\geq \ell$ and hence concludes the proof of parts $1^0$ and $2^0$ of Theorem~\ref{FDH}.

     \medskip
    
  {\bf 4.5.} 
It remains to prove part $3^0$   of Theorem~\ref{FDH}.   We use essentially the same construction of trial functions as in 
part $2^0$ of Theorem~\ref{GG}. Actually, it is slightly more convenient to work with trial functions
  \begin{equation}
 w_{\varepsilon} (\lambda; A)= (  \varepsilon \lambda)^{-1/2} e^{-\varepsilon^{-2} \ln^2 (  \lambda/A) }  .
 \label{eq:GFs}\end{equation}
 If $ A=e^{-a}$, they are
 linked by relation \e{eq:sign1} to functions $u_{\varepsilon} (x; a)$   defined by formula \e{eq:GF} and hence belong to ${\Bbb H}^2_{r}$. The proof below  is significantly more complicated than that of part $2^0$ of Theorem~\ref{GG} because the parameter $A$ in  definition \e{eq:GFs}   will be chosen in rather a special way.
We need the following auxiliary result.

  \begin{proposition}\label{GGF}
 Let the sigma-function $\sigma_{v}$ and trial functions  $ w_{\varepsilon} (\lambda; A)  $ be defined by formulas \e{eq:bvr7} and   \e{eq:GFs}, respectively.
 Then, for any $A >\beta$, we have the relation
   \begin{equation}
\lim_{\varepsilon\to 0} \sigma_{v} [w_{\varepsilon} (A), w_{\varepsilon} (A)]= v_{0}\Gamma(-k)^{-1}\sqrt{\pi/2}
( A - \beta)^{-k-1} e^{-\rho (A -\beta)}.
\label{eq:GGF}\end{equation} 
Moreover, if $B>\beta$, $B\neq A$, then
 \begin{equation}
\lim_{\varepsilon\to 0} \sigma_{v} [w_{\varepsilon} (A), w_{\varepsilon} (B)]= 0.
\label{eq:GGFx}\end{equation} 
 \end{proposition}

 We emphasize that limits \e{eq:GGF} and \e{eq:GGFx} as well as all limits below are uniform with respect to $A$ and $B$ in compact subintervals of $(\beta,\infty)$.
The proof of Proposition~\ref{GGF}
  will be split in several simple lemmas.
  In view of  formula \e{eq:bvr7} and definition \e{eq:di}  we  have
  \begin{equation}
 \sigma_{v} [w_{\varepsilon} (A), w_{\varepsilon} (B)] = \frac{v_{0}e^{\beta\rho}}{\Gamma(-k)}
    \int_{ \beta}^\infty 
  (\lambda - \beta)^{-k-1} \big( \psi_{\varepsilon} (\lambda;A, B)-  \Psi_{\varepsilon}^{(n)} (\lambda ; A,B)   \big) d \lambda
\label{eq:GGF1}\end{equation} 
where
 \begin{equation}
   \psi_{\varepsilon} (\lambda; A,B)=  \varepsilon^{-1} e^{-\rho\lambda} \lambda^{-1} \exp\big( -\varepsilon^{-2} (\ln^2 ( \lambda/A)+ \ln^2 ( \lambda/B))\big),
\label{eq:GGF2}\end{equation}
 \begin{equation}
 \Psi_{\varepsilon}^{(n)} (\lambda; A,B)=   \sum_{p=0}^n \frac{1}{p!} \psi_{\varepsilon}^{(p)} (\beta ; A,B)(\lambda -\beta)^p , \q n=[k]  .
\label{eq:GGF2x}\end{equation}

Our study of integral \e{eq:GGF1}  relies on the following arguments. First, term \e{eq:GGF2x} is important in a neighborhood of the point $\lambda=\beta$ only, and it can be neglected  away from this point. Second, if $A=B$, then  the asymptotics of integral \e{eq:GGF1} as $\varepsilon\to 0$ is determined by a neighborhood of the point $\lambda_{0}= A= B$ (the exponential term in \e{eq:GGF2} equals $1$ at $\lambda_{0}$), but there is no such point if $A\neq B$.

Differentiating definition \e{eq:GGF2}, we obtain bounds on derivatives of function \e{eq:GGF2}.  
 
\begin{lemma}\label{GGF1}
For all $p \in {\Bbb Z}_{+}$ and $\lambda\geq \beta >0$,
  \begin{equation}
|\psi_{\varepsilon}^{(p)} (\lambda; A,B) |\leq C_p (A,B) \varepsilon^{-1-2j} \lambda^{-1}\exp\big( -\varepsilon^{-2} (\ln^2 (  \lambda/ A)+
\ln^2 ( \lambda/B))\big) .
\label{eq:GGF5}\end{equation} 
 \end{lemma} 
 
We suppose for definiteness that $B\leq A$ and put $\lambda_{0} = (A+B)/2$. First, we consider a neighborhood of the point $\lambda=\beta$.

\begin{lemma}\label{GGF3}
Let $\beta< \lambda_{1} < \lambda_{0} $. Then 
 \begin{equation}
  \lim_{\varepsilon\to 0}
    \int_{ \beta}^{\lambda_{1} }
  (\lambda-\beta  )^{-k-1} \big| \psi_{\varepsilon} (\lambda; A,B)-  \Psi_{\varepsilon}^{(n)} (\lambda ; A,B)   \big| d \lambda=0.
\label{eq:GGF6}\end{equation}
 \end{lemma} 
 
 \begin{pf}
Since
  \[
  \big| \psi_{\varepsilon} (\lambda; A,B)-  \Psi_{\varepsilon}^{(n)} (\lambda ; A,B)   \big| \leq (n+1)!^{-1}(\lambda-\beta)^{n+1}\max_{\lambda\in [\beta ,\lambda_{1}]}  |\psi_{\varepsilon}^{(n+1)} (\lambda ; A,B ) |
\]
and $\ln^2 (\lambda_{1}/\lambda_{0})\leq \ln^2 (\lambda / A)$, 
it follows from \e{eq:GGF5}  that
 \[
  \big| \psi_{\varepsilon} (\lambda; A,B)-  \Psi_{\varepsilon}^{(n)} (\lambda ; A,B)   \big| \leq C_{n} (A,B)
  (\lambda-\beta)^{n+1} \varepsilon^{-3-2n} e^{-\varepsilon^{-2} \ln^2 (\lambda_{1}/\lambda_{0})}.
  \]
 This implies \e{eq:GGF6} because $n> k-1$.
\end{pf}

Away from the point $\lambda=\beta$ term \e{eq:GGF2x} is  negligible. 
 
 \begin{lemma}\label{GGF1p}
If $ \lambda_{1} > \beta$, then
\begin{equation}
\lim_{\varepsilon\to 0}    
    \int_{\lambda_{1}}^\infty
  (\lambda - \beta)^{-k-1}  | \Psi_{\varepsilon}^{(n)} (\lambda; A,B)  |d\lambda = 0.
\label{eq:GGF4w}\end{equation} 
 \end{lemma} 
 
  \begin{pf}
  In view of \e{eq:GGF5} where $\lambda=\beta$ the integral here is bounded by
 \[
 C_{n}(A,B)   \exp\big( -\varepsilon^{-2} \ln^2 (  \beta/A) \big) \sum_{p=0}^n       \varepsilon^{-1-2p}  \int_{\lambda_{1}}^\infty
  (\lambda - \beta)^{j-k-1}   d\lambda .
\]
The integrals here are convergent because $p-k\leq n-k <0$, and hence this expression  tends to zero as $\varepsilon\to 0$ because $  \beta<A$.
\end{pf}

 In the next result, we have to distinguish the cases
   $A\neq B$ and $A=B$.

\begin{lemma}\label{GGF2}
 If $A\neq B$, then for any $  \lambda_{1} > \beta$
  \begin{equation}
\lim_{\varepsilon\to 0}    
    \int_{\lambda_{1}}^\infty
  (\lambda - \beta)^{-k-1}   \psi_{\varepsilon} (\lambda; A,B)  d\lambda =  0.
\label{eq:GGF4n}\end{equation} 
 If $A = B$, then for all $\lambda_{1} \in (\beta ,\lambda_{0})$
   \begin{equation}
\lim_{\varepsilon\to 0}    
    \int_{\lambda_{1}}^\infty
  (\lambda - \beta)^{-k-1}   \psi_{\varepsilon} (\lambda; A,A)  d\lambda =  \sqrt{\pi/2}
(A- \beta)^{-k-1} e^{-\rho A } .
\label{eq:GGF4}\end{equation} 
 \end{lemma} 
 
 \begin{pf}
 If $B<A$, then it follows from  \e{eq:GGF5}  that the integral in \e{eq:GGF4n} does not exceed
  \begin{equation}
 C (\lambda_{1}) \varepsilon^{-1 }  \Big(\int_{\lambda_{1}}^{\lambda_0}  \exp\big( -\varepsilon^{-2}  
\ln^2 (  \lambda/ A)\big) d\lambda+ \int_{\lambda_0}^\infty  \exp\big( -\varepsilon^{-2}  
\ln^2 (  \lambda/B)\big) \lambda^{-1} d\lambda\Big). 
\label{eq:GGF4m}\end{equation}
If $\lambda\leq \lambda_{0}$, then $  \lambda < A$ and $-\ln^2 (  \lambda/ A)\leq -\ln^2 (  \lambda_{0}/ A)$. Therefore the first term in \e{eq:GGF4m} is estimated by $C  \varepsilon^{-1 }  \exp\big( -\varepsilon^{-2}  
\ln^2 (  \lambda_{0}/ A)\big)$ which tends to zero as $\varepsilon\to 0$ because $  \lambda_{0} < A$.

  Making the change of variables $\lambda= B e^{\varepsilon x }$, we see that the second term in
 \e{eq:GGF4m}    is estimated by the integral of   $e^{-x^2}$ over the interval $(x_{0}(\varepsilon),\infty)$
 where $x_0(\varepsilon)= \varepsilon^{-1}\ln(\lambda_0/ B)\to +\infty$ as $\varepsilon\to 0$  because $  \lambda_{0} > B$.

% In the case $a=b$,  the asymptotics of the integral in \e{eq:GGF4} is determined   by a neighborhood of the point $ \lambda_{0}$.
   
    Making the   change of variables $\lambda= A  e^{\varepsilon x  }$, we see that   integral  \e{eq:GGF4}  equals
   \[ 
 \varepsilon^{-1}   \int_{\lambda_{1}}^\infty
  (\lambda - \beta)^{-k-1}      e^{-\rho\lambda} e^{ -2\varepsilon^{-2} (\ln^2 ( \lambda/A)}   \lambda^{-1} d\lambda = 
    \int_{x_{1}(\varepsilon)}^\infty
  ( Ae^{ \varepsilon x} - \beta)^{-k-1}     e^{-\rho Ae^{  \varepsilon x}}  e^{-2x^2} dx
  \]
where $x_1(\varepsilon)=\varepsilon^{-1}\ln (\lambda_1 / A)\to -\infty$   as $\varepsilon\to 0$ because $  \lambda_1 < A$.
Therefore the right-hand side here converges as
 $\varepsilon\to 0$ to the corresponding integral over $\Bbb R$. It equals the right-hand side of \e{eq:GGF4}.
\end{pf}

% It remains to consider a neighborhood of the point $\lambda= \beta$.

Let us return to representation \e{eq:GGF1}. In view of Lemma~\ref{GGF3} the integral over $(\beta, \lambda_{1})$ can be neglected for any $\lambda_{1}<  (A+B)/2$. If $A\neq B$, then the integral \e{eq:GGF1} over $( \lambda_{1}, \infty)$ also tends to zero according to relations \e{eq:GGF4w} and \e{eq:GGF4n}. This yields \e{eq:GGFx}. In the case $A=B$ the integral \e{eq:GGF1} over $( \lambda_{1}, \infty)$  tends to the right-hand side of \e{eq:GGF4} according to relations \e{eq:GGF4w} and \e{eq:GGF4}.
This yields \e{eq:GGF}
and concludes the proof of   Proposition~\ref{GGF}.  

  \medskip
    
  {\bf 4.6.} 
Now we are in a  position to prove part $3^0$   of Theorem~\ref{FDH}. Let the functions $w_{\varepsilon} (\lambda; A ) $ be defined by formula \e{eq:GFs}. Observe that  under Assumption~\ref{Hfree} we have the estimate
\[
\sigma_{0} [w_{\varepsilon}(A) , w_{\varepsilon} (A)] =\varepsilon^{-1}\int_{-\infty}^\infty s_{0} (x) e^{-2 \varepsilon^{-2}(x-a)^2} dx \leq C \varepsilon^{-1}\int_{-\infty}^\infty  e^{(l+ |a|) |x|} e^{-2 \varepsilon^{-2}x^2} dx
\]
where $a=-\ln A$ and $l=\max\{ l_{-}, l_{+}\}$. Making here the change of variables $x=\varepsilon y$, we see that this expression is bounded by a constant $c_{0}>0$ which does not depend on $\varepsilon\in (0,1]$ and on the parameter $a $ in  a compact interval of $\Bbb R$.

Let $N$ be given. We look for trial functions in the form
 \[
w_{\varepsilon} (\lambda)=\sum_{j=1}^N \nu_{j} w_{\varepsilon} (\lambda; A_{j})
\]
where  $A_{j}  > \beta$ for all $j=1,\ldots, N$ and $\nu_1, \ldots, \nu_{N}$ are arbitrary complex numbers. As we have seen, 
  \begin{equation}
\sigma_{0} [w_{\varepsilon} , w_{\varepsilon} ] \leq  N \sum_{j=1}^N |\nu_{j} |^2 \sigma_{0} [w_{\varepsilon}(A_{j}) , w_{\varepsilon} (A_{j})] \leq c_{0}N \sum_{j=1}^N |\nu_{j} |^2.
\label{eq:TF1}\end{equation}

Next, we consider the form $\sigma_{v} [w_{\varepsilon} , w_{\varepsilon} ] $. Observe that under our assumptions $v_{0}\Gamma (-k)< 0$.   We choose the points $A_{1},\ldots, A_{N}$  so close to $  \beta$ that
\[
v_{0} \Gamma(-k)^{-1}\sqrt{\pi/2}
(A_{j}- \beta)^{-k-1} e^{-\rho (A_{j} -\beta)}\leq -3 c_{0}, \q j=1,\ldots, N.
\]
Then it follows from Proposition~\ref{GGF} that for a sufficiently small $\varepsilon>0$
\begin{equation}
\sigma_{v} [w_{\varepsilon} , w_{\varepsilon} ]  \leq -2 c_{0}\sum_{j=1}^N |\nu_{j} |^2.
\label{eq:TF3}\end{equation}

Comparing estimates \e{eq:TF1} and \e{eq:TF3}, we see that
\[
\sigma  [w_{\varepsilon} , w_{\varepsilon} ]   \leq -  c_{0}\sum_{j=1}^N |\nu_{j} |^2.
\]
for arbitrary numbers $\nu_1,\ldots, \nu_{N}\in{\Bbb C}$. In view of formula \e{eq:DHW} where $w_{\varepsilon}$  and $f_{\varepsilon}$ are related by  equation \e{eq:sign2}, this yields us the subspace in ${\cal D} [h]$ of dimension $N$ where the form $h[f_{\varepsilon},f_{\varepsilon}] $ is negative. This concludes the proof of   Theorem~\ref{FDH}.

  %************************************************************
  %************************************************************
\section{Quasi-Carleman and finite rank Hankel operators}  
 %*****
  
  Here we consider perturbations of singular quasi-Carleman operators $H_{0}$ by finite rank self-adjoint Hankel operators $V$. 
  We shall prove that $N_{-} (H_{0} +V )=N_{-} (V)$, that is,  adding $H_{0}$ to $V$ does not change
    the total number of   negative   eigenvalues of a finite rank Hankel  operator $V$. Our proof here is relatively similar to 
  that of  Theorem~\ref{FDH}, but new difficulties   arise because the singularities of the sign-function $s_{v} (x)$ may lie in the complex plane; in this case $s_{v} (x)$ is even more singular than function \e{eq:bvr7}. 
  
  % In contrast to Section~4 it is now slightly more convenient to work with  sign-functions $s(x)$ instead of sigma-functions $\sigma(\lambda)$.

\medskip

{\bf 5.1.}
The unperturbed operator $H_{0}$ is the same as in Section~4. Thus we accept  
Assumption~\ref{Hfree} and  define  $H_{0}$ by its quadratic form
\e{eq:CB}.

  Recall that integral kernels of finite rank Hankel operators $V$  are given (this is the classical  Kronecker theorem -- see, e.g., Sections~1.3 and 1.8 of the book \cite{Pe}) by the formula
\begin{equation}
v(t) =\sum_{m=1}^M  P_{m} (t)e^{- \beta_{m}  t}
\label{eq:FDvm}\end{equation}
where $\Re \beta_{m}  >0$ and 
 $   P_{m} (t)$ 
    are polynomials. We consider   self-adjoint $V$ when $ v(t)=\ov{v(t)} $. If  $\Im \beta_{m}  \neq 0$,  then necessarily the sum in  \e{eq:FDvm} contains both terms  $P_{m} (t) e^{- \beta_{m}  t}$ and   $\ov{P_{m} (t)} e^{- \bar{\beta}_{m}  t}$. Let $\Im \beta_{m}=0$ for $m=1,\ldots, M_{0}$   and $\Im \beta_{m}>0$,  $\beta_{ M_{1}+m}=\bar{\beta}_{ m}$ for $m=M_{0} +1, \ldots,M_{0} + M_{1} $. Thus $M=M_{0}+2M_{1}$; of course the cases $M_{0}=0$ or $M_{1}=0$ are not excluded.  We have $ P_{m} (t) = \ov{P_{m} (t)} $ for $m=1,\ldots, M_{0}$ and $ P_{M_{1}+m} (t) = \ov{P_{m} (t)}$  for $m=M_{0} +1, \ldots,M_{0} + M_{1} $.  Let $K_{m}=\deg P_{m}$. Then   
$
\rank V =\sum_{m=1}^M K_{m}+M.
$

     For $m=1,\ldots, M_{0}$, we
 set
\[
 {\sf p}_{m}=P_{m}^{(K_m)},
\]
that is, ${\sf p}_{m}/K_{m}!$ is the coefficient at $t^{K_{m}}$ in the polynomial  $P_{m} (t)$, 
 and
    \begin{equation}
\left. \begin{aligned}
  {\cal N}_m &=(K_{m} +1)/2 \; &&{\rm if} \;  K_{m} \; {\rm is  } \; {\rm   odd}
 \\
 {\cal N}_m &= K_{m}  /2 \; &&{\rm if} \; K_{m} \; {\rm is  } \; {\rm   even} \; {\rm   and}
 \q {\sf p}_{m}  > 0
 \\
 {\cal N}_m &= K_{m}  /2  +1\; &&{\rm if} \; K_{m} \; {\rm is  } \; {\rm   even} \; {\rm   and}
 \q {\sf p}_{m} < 0 .
    \end{aligned}
    \right\}
    \label{eq:RXm}\end{equation}

   Our main result is formulated as follows.
 
  \begin{theorem}\label{FDH1}
  Let a function $s_{0}(x)$ satisfy Assumption~\ref{Hfree}, and let $H_{0}$ be the self-adjoint positive operator defined by the quadratic form \e{eq:CB}.
 Let $V$ be the self-adjoint Hankel operator  of  finite rank with kernel $v(t)$ given by formula \e{eq:FDvm}, and let the numbers ${\cal N}_m$ be defined by formula \e{eq:RXm}.  
  Then the total number $N_{-} (H)$  of   $($strictly$)$   negative eigenvalues of the operator $ H=H_{0}+V$ is given by the formula 
   \begin{equation}
 N_- (H) =  \sum_{m=1}^{M_{0}}  {\cal N}_m+ \sum_{m=M_{0}+1}^{M_{0}+M_1}  K_{m}  +M_{1} =:{\cal N} .
   \label{eq:TN}\end{equation}
 \end{theorem}
 
 Theorem~\ref{FDH1} generalizes the corresponding result of \cite{Yf} where the sign-function $s_{0}(x)$ was supposed to be bounded; in this case the operator $H_{0}$ is also bounded. In particular, formula \e{eq:TN} was established in \cite{Yf} in the case $H_{0}=0$, that is, for finite rank Hankel operators $H$.
  We emphasize that the right-hand side of \e{eq:TN} does not depend on the operator $H_{0}$.  Therefore using \e{eq:TN} for $H_{0}=0$, we obtain equality \e{eq:FDb}.

 Since $H_{0}\geq 0$, we have $N_{-} (H) \leq N_{-} (V)={\cal N} $. Thus we only have   to prove that
  \begin{equation}
 N_{-} (H) \geq {\cal N}  .
   \label{eq:TNC}\end{equation}
  To that end, we construct trial functions $f$ such that $h[f,f]< 0$. We emphasize that the constructions for the terms in \e{eq:FDvm}  corresponding to $\Im \beta_{m}=0$ and to $\Im \beta_{m}\neq 0$ are essentially different.

  \medskip
    
  {\bf 5.2.}
  In  this subsection we collect some results of \cite{Yf} which we use below.
First we recall the explicit  expression  for the sign-function $s_{v} (x)$ of the Hankel operator with  kernel $v(t) =t^j e^{-\beta t}$. 

 \begin{lemma}\label{Fr}
     Let $v(t) =t^j e^{-\beta t}$  where $\Re \beta  >0$ and  $j \in {\Bbb Z}_{+}$. 
      If $j=0$, then  
        \begin{equation}
s_{v} (x) = \beta^{-1}\delta (x -\kappa), \q \kappa=-\ln \beta, \q -\pi/2 < \Im \kappa <\pi/2.
\label{eq:E6}\end{equation}
 If $j\geq 1$, then
  \begin{equation}
s_{v} (x) =  \beta^{-1-j}(1 -\partial)\cdots (j -\partial)\delta (x -\kappa).
\label{eq:E7}\end{equation}
 \end{lemma} 
 
  Clearly, $s_{v}\not\in {\cal S}'$   unless $\Im \beta
 =0$, but the corresponding sigma-function $\sigma_{v}\in {\cal Y}'$. Actually,  distributions \e{eq:E6} and \e{eq:E7} are well defined as antilinear functionals on test functions $u(z)$  analytic in the strip $-\pi/2 < \Im z < \pi/2$. We put $u^* (z)=\ov{u(\bar{z})}$. It follows from Lemma~\ref{Fr} that $s_{v}[u,u]=\la s_{v }, u^* u\ra$ is determined by values of $u(z)$ and its derivatives at the points $z=\kappa$ and $z=\bar{\kappa} $. To be more precise, Lemma~\ref{Fr} implies the following 
 result.

 \begin{lemma}\label{FR}
 Let $ v(t) =  P (t)e^{- \beta  t} $
where $\Re \beta >0$ and $P (t)= p_{K} t^K +\cdots$, $p_{K}  \neq 0 $, is a polynomial of degree $K$.
Put\footnote{The upper index ``$\top$" means that a vector is regarded as a column.}
      \begin{equation}
J_{K} (\kappa) u= (u (\kappa  ),u' (\kappa  ), \ldots, u^{(K )}(\kappa ))^\top \in {\Bbb C}^{K +1}.
 \label{eq:Z}\end{equation}
      Then 
     \[
  s_{v} [u, u]  =  ({\bf S} (P, \kappa)  J_{K} (\kappa) u ,   J_{K} (\bar\kappa) u )_{K+1} ,\q \kappa=-\ln \beta,  
 \]
  where   $( \cdot,  \cdot )_{K+1} $ is   the scalar product  in   ${\Bbb C}^{K +1}$ and the matrix ${\bf S}(P, \beta)$ is skew triangular,
  that is, its entries $s_{j, \ell}    =0$    for $j+\ell   > K $. Moreover, we have
   \begin{equation}
   s_{j, \ell}    = C_{K }^j \beta^{-1-K}  p_{K } \q {\rm for}\q j+\ell   = K 
 \label{eq:E81vvx}\end{equation}
 and
${\bf S} (\bar{P},\bar\beta)= {\bf S} (P,\beta)^* $.
 \end{lemma}

 According to \e{eq:E81vvx} we have $\det{\bf S}\neq 0$. 
  Actually,  all entries $s_{j, \ell}    $    of the matrix $ {\bf S} (P, \beta)$  (we call it the {\it sign-matrix}  of the kernel  $v(t)$) admit simple expressions in terms of the coefficients of the polynomial $P(t)$, but below we need only the information collected in Lemma~\ref{FR}.
  
  In the symmetric case,  we use       the following result on spectra of the sign-matrices.  
      
       \begin{lemma}\label{ST}
Let $\beta=\bar{\beta}$ and $P(t)=\ov{P(t)}$.   If $K$ is odd, then $ {\bf S} (P,\beta) $ has $(K+1)/2$ positive and $(K+1) /2$ negative eigenvalues. If $K$ is even, then $ {\bf S} (P,\beta)$ has $K/2+1$ positive and $K/2$ negative eigenvalues for $ {\sf p}_{K }>0$ and it has $K/2 $ positive and $K/2+1$ negative eigenvalues for $ {\sf p}_{K }<0$.
 \end{lemma}

       In the complex case, Lemma~\ref{FR} implies the following assertion.

  \begin{proposition}\label{FR1qc}
  Let 
   \[
v(t) =  P (t)e^{- \beta  t} +  \ov{P (t)}e^{- \bar\beta  t}, \q \Re\beta>0, \q \Im\beta>0.
\]
Put
  \begin{equation}
\wt{J}_{K} (\kappa) u =  (  J_{K} (\kappa) u , J_{K} (\bar\kappa) u )^\top.
\label{eq:FDvct}\end{equation}
Then
   \[
  s_{v}[u,u]   = (  {\bf \wt{S}} (P, \beta) 
 (  \wt{J}_{K} (\kappa) u ,
  (  \wt{J}_{K} (\kappa) u )_{2K+2}   ,    \q \kappa=-\ln \beta, 
\]
       where the sign-matrix
          \begin{equation}
   {\bf \wt{S}} (P, \beta) = \begin{pmatrix}
0 & {\bf S} (P,\beta)^*
\\
{\bf S} (P,\beta)    & 0
\end{pmatrix}.    
 \label{eq:sicc1}\end{equation}
        \end{proposition}

        Obviously, the spectrum of matrix \e{eq:sicc1} is symmetric so that it consists of $K+1$ positive and $K+1$ negative
        eigenvalues.
        
        Let us   collect the results obtained together.
        
          \begin{theorem}\label{FINR}
          Let $v(t)$ be kernel  \e{eq:FDvm}, let the operators $J_{K_{m}} (\kappa_{m})$ and $\wt{J}_{K_{m}} (\kappa_{m})$ be defined by formulas \e{eq:Z} and \e{eq:FDvct}, respectively, and let ${\bf S} (P_{m}, \beta_{m}) $ and 
          ${\bf \wt{S}}  (P_{m}, \beta_{m}) $ be the corresponding sign-matrices. 
   Then for all functions $u(z)$ analytic in the strip $-\pi/2< \Im z< \pi/2$, we have
          \begin{align}
 s_{v}[u , u] = & \sum_{m=1}^{M_{0} } ({\bf S} (P_{m}, \beta_{m})  J_{K_{m}} (\kappa_{m}) u ,   J_{K_{m}} ( \kappa_{m}) u )_{K_{m}+1}
\nonumber \\
& + \sum_{m=M_{0}+ 1}^{M_{0}+M_{1}}  (  {\bf \wt{S}} (P_{m}, \beta_{m}) 
 \wt{J}_{K_{m}} (\kappa_{m}) u ,    
  \wt{J}_{K_{m}} (\kappa_{m}) )_{2K_{m}+2} .
 \label{eq:si2}\end{align}
           \end{theorem}

  \medskip
    
  {\bf 5.3.}
  For the construction of trial functions $u(x)$ where the form $s[ u,u]<0$,   we need the following assertion. We emphasize that the considerations of real and complex $\kappa_{m}$ in \e{eq:si2} are essentially different.
  
  %Let us now calculate the  form $s [u,u]$ on some special functions $u(x)$ of gaussian type. 

 \begin{lemma}\label{zzb}
  Let $\kappa \in {\Bbb C}$, $\varepsilon>0$ $(\varepsilon$ is a small parameter$)$ and let $\omega(z)$ be a polynomial such that $\omega(\kappa)\neq 0$. If $\kappa=\bar{\kappa}$, we set 
    \begin{equation}
  \varphi (z;\varepsilon)=\omega (z)   e^{-(z - \kappa )^2/\varepsilon^2}.
\label{eq:AK1}\end{equation}
 If $\kappa= \kappa' + i \kappa''$ where $\kappa'' \neq 0$, we set 
    \begin{equation}
  \varphi (z;\varepsilon)=\omega (z) e^{-i \sgn\kappa'' (z-\kappa)/\varepsilon}  e^{-(z - \kappa' )^2 }.
\label{eq:AK2}\end{equation}
Let $a_{0},a_{1},\ldots, a_{K}$ be any given numbers. Then there exists a polynomial
 \begin{equation}
 Q  (z ;\varepsilon )= \sum_{p=0}^{K } q_p (\varepsilon)(z - \kappa )^p
\label{eq:AF4}\end{equation}
such that the function
 \begin{equation}
 \psi (z;\varepsilon)= Q (z  ;\varepsilon)\varphi (z;\varepsilon) 
\label{eq:AF2b}\end{equation}
satisfies the conditions
 \begin{equation}
\psi ^{(j)}(\kappa ; \varepsilon) =a_{j}, \q j =0,1, \ldots, K.
\label{eq:AK3}\end{equation} 
Moreover, the coefficients of polynomial \e{eq:AF4} satisfy estimates  
\begin{equation}
|q_p ( \varepsilon)|Ê\leq C  \varepsilon^{-p},\q p =0,1, \ldots, K.
\label{eq:AG5b}\end{equation}
 \end{lemma}
  
 \begin{pf}
 In view of  \e{eq:AF4}, \e{eq:AF2b}  conditions \e{eq:AK3} yield the equations 
\begin{equation}
 j!  \sum_{p=0}^j    (p-j)!^{-1}  q_p (\varepsilon)\varphi ^{(j-p)}(\kappa ;\varepsilon) =a_{j}, \q j=0, 1,\ldots, K , 
\label{eq:AF5}\end{equation}
for the coefficients $q_p (\varepsilon)$. Let us consider these equations successively starting from $j=0$.
Observe that $\varphi(\kappa ; \varepsilon) = \varphi (\kappa ; 0) =  \omega(\kappa)e^{\kappa''^2}\neq 0$. Therefore
\[
q_{0 } = \varphi (\kappa ; 0) ^{-1} a_{0}  .
\]
Then equation \e{eq:AF5} determines $q_{j } (\varepsilon)$ if $q_{0 }, q_{1 }(\varepsilon), \ldots, q_{j-1 }(\varepsilon)$ are already found: 
\[
 \varphi (\kappa ; 0) q_j (\varepsilon)=  j!^{-1}a_{j}-
  \sum_{p=0}^{j-1}   (p-j)!^{-1} q_p (\varepsilon)\varphi ^{(j-p)}(\kappa ;\varepsilon)  .
\]
Since for both functions \e{eq:AK1} and \e{eq:AK2}
 \[
|\varphi ^{(k)}(\kappa ; \varepsilon) |\leq C  \varepsilon^{-k},  
\]
estimates \e{eq:AG5b} on $q_0 ( \varepsilon), \ldots, q_{j-1} ( \varepsilon)$ imply the same estimate on $q_j ( \varepsilon)$.   
 \end{pf}
 
 If $z=x\in{\Bbb R}$, then functions  \e{eq:AK1} and \e{eq:AK2}    satisfy estimates
  \begin{equation}
 | \varphi (x;\varepsilon)|Ê\leq C (1+ |x-\kappa|^{\deg \omega})    e^{-(x - \kappa )^2/\varepsilon^2 },
 \q \kappa=\bar{\kappa},
 \label{eq:BL1}\end{equation}
and
 \begin{equation}
 | \varphi(x;\varepsilon)|Ê\leq C e^{-| \kappa''|/\varepsilon }  (1+ |x-\kappa'|^{\deg \omega})    e^{-(x - \kappa' )^2 },
 \q \kappa''\neq 0,
 \label{eq:BL2}\end{equation}
 respectively.  Here we have taken into account that
  \[
  |e^{-i \sgn\kappa'' (z-\kappa)/\varepsilon} |=  e^{-| \kappa''|/\varepsilon } .
  \] 
  Note also that in view of \e{eq:AG5b} polynomial \e{eq:AF4}  satisfies the estimate
   \begin{equation}
| Q(x;\varepsilon)|Ê\leq C (1+ |x-\kappa'|^K \varepsilon^{-K}).
  \label{eq:BLL}\end{equation}
This leads to the following assertion.

\begin{corollary}\label{zzbc}
Under the assumptions of Lemma~\ref{zzb}, for any $l\in {\Bbb R}$, 
 \begin{equation}
 \int_{-\infty}^\infty e^{l |x|} |\psi (x; \varepsilon)|^2 dx = O (\varepsilon),\q \varepsilon\to 0 .
 \label{eq:YY6}\end{equation}
  \end{corollary}

\begin{pf}
If $\kappa'' =0$, we use estimates \e{eq:BL1} and \e{eq:BLL}.
  Making the change of variables $x-\kappa =\varepsilon y$, we see that
    integral \e{eq:YY6} is bounded by $C \varepsilon $. If $\kappa''\neq 0$, then this integral tends to zero exponentially due to the factor $e^{-| \kappa''|/\varepsilon } $ in  \e{eq:BL2}.
\end{pf}

   \medskip
    
  {\bf 5.4.}
  Let us return to Hankel operators $H=H_{0}+V$. For kernels \e{eq:FDvm}, we put $\kappa_{m}= -\ln \beta_{n}$. Then $\kappa_{m}=\bar{\kappa}_{m}$ for $m=1,\ldots, M_{0}$ and $\kappa_{m}=\bar{\kappa}_{m+M_{1}}$ for $m=M_{0}+1,\ldots, M_{0}+M_{1}$. For all $k=0,\ldots, K_m$,  $m=1,\ldots, M $,
     we construct the functions $\psi_{k, m}(z,\varepsilon)$  by formulas of Lemma~\ref{zzb} where $\kappa=\kappa_{m}$.    Of course we use formula  \e{eq:AK1} if $\Im \kappa_{m}=0$ and \e{eq:AK2} if $\Im \kappa_{m}\neq 0$.
     We require that $\psi_{k, m}^{(l)}(\kappa_m; \varepsilon) = \d_{k, l}$ for all  $k,l=0,\ldots, K_m$.  Let us set
   $\omega_{1} (z)=1$ if $M=1$ and
   \[
\omega_{m} (z)=   \prod_{n=1; n\neq m}^M (z- \kappa_{n})^{K_{n} +1}\q {\rm if }\q M\geq 2.
\]
Then $\omega_{m} (\kappa_{m}) \neq 0$ and due to this factor in \e{eq:AK1} and \e{eq:AK2}  the function $\psi_{k, m} (z; \varepsilon)  $ satisfies   the conditions 
$\psi_{k, m}^{(l)}(\kappa_n; \varepsilon) = 0$ for all  $n\neq m$ and $l=0,\ldots, K_n$.  By virtue of the conditions at the points
$\kappa_{m}$, $m=1,\ldots, M $,   for every fixed $\varepsilon>0$   all  functions $\psi_{k, m} (z;\varepsilon)$  
    are linearly independent.
    
     For arbitrary complex numbers $\nu_{k,m}$, we put
\begin{equation}
u_{m}(z;\varepsilon)=   \sum_{k=0}^{K_m} \nu_{k,m}\psi_{k, m} (z;\varepsilon)  
\label{eq:YY}\end{equation}
for $   m=1,\ldots, M_{0}$
and
\begin{equation}
u_{m}(z;\varepsilon)=   \sum_{k=0}^{K_m} \big( \nu_{k,m}\psi_{k, m} (z;\varepsilon) + \nu_{k,m+M_{1}}\psi_{k, m+M_{1}} (z;\varepsilon)  \big) 
\label{eq:YYY}\end{equation}
for $  m=M_{0}+1,\ldots, M_{0}+M_{1}$. The  functions $u_{1} (z;\varepsilon), \ldots, u_{M_{0}+M_{1}} (z;\varepsilon)$  
    are of course linearly independent.

    Let ${\sf a}_{m}= (\nu_{0,m}, \nu_{1,m}, \ldots, \nu_{K_{m},m})^\top \in {\Bbb C}^{K_{m}+1}$ for $m= 1,\ldots, M $. We put
   ${\bf a}_{m}={\sf a}_{m}$,  ${\bf J}_{m} u= J_{K_{m}} (\kappa_{m}) u$, $ {\bf  S}_{m} = {\bf  S} (P_{m}, \beta_{m}) $  for $m= 1,\ldots, M_{0} $ and  ${\bf a}_{m}=({\sf a}_{m}, {\sf a}_{m+M_{1}})^\top$,  ${\bf J}_{m} u= \wt{J}_{K_{m}} (\kappa_{m}) u$, ${\bf  S}_{m} =  {\bf \wt{S}} (P_{m}, \beta_{m}) $   for $m= M_{0}+ 1,\ldots, M_{0} +M_{1}$. Note that each matrix ${\bf  S}_{m}$ has ${\cal N}_{m}$ negative eigenvalues with ${\cal N}_{m}$ defined by  formula \e{eq:RXm}.
   For  $  m= 1,\ldots, M_{0} $, this result follows from  Lemma~\ref{ST}.  For  $  m= M_{0}+1,\ldots, M_{0} +M_1$, this is a direct consequence of representation \e{eq:sicc1}.

For functions  \e{eq:YY} and  \e{eq:YYY}, we  have ${\bf J}_{m}u_{m}(\varepsilon)= {\bf a}_{m}$ and ${\bf J}_{n}u_{m}(\varepsilon)=0$ if $n\neq m$. Therefore it follows from formula  \e{eq:si2}  that 
  \[
 s_{v} [u_{m} (\varepsilon) , u_{m} (\varepsilon)] =  ( {\bf S}_{m} {\bf a}_{m}   ,  {\bf a}_{m} ), \q  s_{v} [u_{m} (\varepsilon) , u_n (\varepsilon)] = 0, \q \forall \varepsilon>0,   
 \]
 for all $  m,n= 1,\ldots, M_{0}+M_{1}$, $n\neq m$.
 According to Corollary~\ref{zzbc} under Assumption~\ref{Hfree} we have
  \[
 s_{0} [u_{m} (\varepsilon) , u_n (\varepsilon)]  \leq  C\varepsilon  \| {\bf a}_{m} \|  \| {\bf a}_n \|, \q \forall
 m,n= 1,\ldots, M_{0}+M_{1} .
        \] 
 Thus if     ${\bf a}_{m} $ is an eigenvector of the matrix ${\bf S}_{m}$ corresponding to its negative eigenvalue $-\mu_{m}$, then  
          \[
 s [u_{m} (\varepsilon) , u_{m} (\varepsilon)]  \leq -(\mu_{m} - C\varepsilon ) \| {\bf a}_{m} \|^2 .
        \]  
        This expression is negative if $C\varepsilon<\mu_{m}$.    This yields the estimate 
        \begin{equation}
         s [u (\varepsilon) , u (\varepsilon)] \leq -c \| u\|^2, \q c>0,
\label{eq:YZY}\end{equation}
         for all  linear combinations $u (z; \varepsilon)$  of functions $u_{m} (z; \varepsilon)$. 
         Since we have ${\cal N}_{m}$ linearly independent functions $u_{m} (z; \varepsilon)$, the dimension of their linear combinations  equals ${\cal N}={\cal N}_{1}+\cdots + {\cal N}_{M_{0}+M_{1}}$.

Let $f (t; \varepsilon) $  be the  functions  linked to   $u (x; \varepsilon)$   by formula  \e{eq:sign2}.  Since the functions $ \psi_{k,m} (x;\varepsilon)$ defined by equalities  \e{eq:AK1}, \e{eq:AK2} and
\e{eq:AF2b} decay super-exponentially as $|x|\to \infty$, the functions
  $f  (\varepsilon)\in L^2 ({\Bbb R}_{+})$. These functions belong to ${\cal D} [h]$ because $ s [u  (\varepsilon) , u  (\varepsilon)] <\infty$.  It now follows from identity \e{eq:gB} and estimate \e{eq:YZY} that
\[
h[f,f]=\| \sqrt{H}_{0} f (\varepsilon)\|^2 +(Vf  (\varepsilon),f_{m} (\varepsilon))= s [u  (\varepsilon) , u  (\varepsilon)] <0
\]
for all $f ( \varepsilon) \neq 0$  
in  the linear space   of dimension ${\cal N} $. 
This yields  estimate \e{eq:TNC} and thus
   concludes the proof of Theorem~\ref{FDH1}.

 \begin{remark} 
 The condition  $s_{0}\in L^\infty_{\rm loc} ({\Bbb R})$ of regularity of the sign-function $s_{0} (x)$  in Assumption~\ref{Hfree} cannot be omitted. Indeed, consider, for example, the kernels $h_{0}(t)= 2 e^{-t}$, $v (t)= - e^{-t}$. Then $N_{-}(V)=1$ while 
$N_{-}(H)=0$.
    \end{remark}

\appendix{}

%%%%%%%%%%%%%%%%%%%%%%%

%%%%%%%%%%%%%%%%%%

\section{Sandwiched Fourier transforms}

%%%%%%%%%%%%%%%%%%%%%%

 For absolutely continuous measures $dM(\lambda)=\sigma (\lambda) d\lambda$, the results of  Section~3 on the operators ${\sf L}: L^2\to L^2 (M)$ can be reformulated  in terms of the operators
  \begin{equation}
{\cal A} =s (x) \Phi^* v(\xi) 
\label{eq:BF}\end{equation}
where $v(\xi) =\Gamma (1/2- i\xi)$ and $s(x)= \sigma (e^{-x})$. Now we study   operators ${\cal A}$ in the space $L^2: = L^2 ({\Bbb R})$ for sufficiently arbitrary functions $v(\xi)$ and  $s(x)$. Properties of operators \e{eq:BF} for $v$ and $\sigma$ belonging to some spaces $L^p  $ were extensively discussed in the literature (see, e.g., the book \cite{RS}). Here we consider operators \e{eq:BF} with  functions $ s(x)$ growing at infinity. Our goal is to define operators 
\e{eq:BF} as closed (unbounded) operators in $L^2$.

In contrast to Section~3,  our assumptions on $s(x)$ exclude its exponential growth at infinity. We now suppose that 
\begin{equation}
| s (x) |\leq C (1+|x|)^K
\label{eq:BF1}\end{equation}
for some $K\in {\Bbb R}$. With respect to $v(\xi)$, we assume that $v\in C^\infty$ and
\begin{equation}
| v^{(n)}(\xi) |\leq C_{n} (1+|\xi|)^{k_{n}}
\label{eq:BF2}\end{equation}
for all $n=0,1,\ldots$ and some numbers $k_{n}\in {\Bbb R}$. Obviously, for $f\in {\cal S}$, we have $vf\in {\cal S}$, $\Phi^* (vf)\in {\cal S}$ and hence $s\Phi^* (vf)\in L^2 $. It means that ${\cal A} : {\cal S}\to L^2$. 
If $g\in L^2$, then $\bar{s} g\in {\cal S}'$, $\Phi (\bar{s} g)\in {\cal S}'$ and 
\begin{equation}
{\cal A}_{*} g : =\bar{v} \Phi \bar{s} g\in {\cal S}'.
\label{eq:SST}\end{equation}
 It means that    ${\cal A}_{*} : L^2 \to {\cal S}'$.  Moreover, for all $ f \in {\cal S} $ and all $g\in L^2  $, we have the identity
 \begin{equation}
 \la{\cal A} f, g  \ra  =    \la f, {\cal A}_{*} g  \ra. 
\label{eq:LTM1B}\end{equation}

Define now the operator $A_{0} $  in $L^2$ on the domain  ${\cal D}(A_{0})=\cal S$ by the equality $A_{0}f= {\cal A}f$. Let us construct its adjoint operator $A_{0}^*   $. Let ${\cal D}_{*}\subset   L^2  $ consist of   $g \in   L^2  $ such that ${\cal A}_{*} g \in L^2$.  
  
     \begin{lemma}\label{LTMB}
     Under assumptions \e{eq:BF1} and \e{eq:BF2} the operator $A_{0}^*$ is given by the equality $A_{0}^* g ={\cal A}_{*}g$ on the domain ${\cal D}(A_{0}^*) ={\cal D}_{*}$. 
    \end{lemma}
    
     \begin{pf}
     If $ f \in {\cal D}(A_{0}) $ and   $g\in {\cal D}_{*} $, then it follows from identity \e{eq:LTM1B} that
      \[
 (A_{0} f, g)  =   (f, {\cal A}_{*} g).
\]
Hence $   g \in {\cal D}(A_{0}^*)  $  and  $A_{0}^* g ={\cal A}_{*}g$.
     
     Conversely, if $g\in {\cal D}(A_{0}^*) $, then $| ({\cal A} f, g)  | \leq C \| f\| $ for all $ f \in {\cal S} $. In view  of 
     \e{eq:LTM1B}, this estimate implies that ${\cal A}_{*} g \in L^2$ and hence ${\cal D}(A_{0}^*)   \subset {\cal D}_{*}$.
          \end{pf}
  
  \begin{corollary}\label{LTM1B}
  Suppose additionally that $k_{0}=0$ in condition \e{eq:BF2}.
Then the operator $A_{0} $   admits the closure.  
  \end{corollary}
    
      \begin{pf}
      Indeed,   for $g\in {\cal S}$, we have $\bar{s} g\in L^2$ so that ${\cal A}_{*} g\in L^2$ if the function $v$ is bounded. It follows that
            ${\cal S}\subset {\cal D}_{*} ={\cal D}(A_{0}^*) $,  and hence the operator $A_{0}^*$ is densely defined. 
                    \end{pf} 
        
 Next, we construct the second adjoint $ A_{0}^{**}$. Observe that ${\cal A} : L^2  \to  {\cal S}'$.  Let the operator $A$ be defined   by the equality $A f={\cal A } f$ on the domain   ${\cal D} (A)$ which  consists of  all $f \in   L^2  $ such that ${\cal A } f\in L^2$. 
        
        \begin{lemma}\label{LTM2B}
  Under assumptions \e{eq:BF1} and \e{eq:BF2} where $k_{0}< -1/2$ , the inclusion  $A_{0}^{**}\subset A $ holds. 
    \end{lemma}
    
     \begin{pf}
      For $ f \in L^2 $ we have $v f \in L^1 $, and for $g\in \cal S$ we have $\bar{s} g \in L^1 $. Therefore 
        according  to the Fubini theorem,  the identity
  \e{eq:LTM1B} is now true
for all $ f \in L^2 $ and all $g\in \cal S$.  If $f \in {\cal D} (A_{0}^{**})$, then $|(f, {\cal A}_{*}g)|\leq C \| g\|$  for all $g\in {\cal D}_{*}$ and, in particular, for $g\in {\cal S}$. Thus it follows from    \e{eq:LTM1B} that $|({\cal A}  f,  g)|\leq C \| g\|$ and hence ${\cal A} f\in L^2  $. Moreover,   $A_{0}^{**}f= {\cal A} f$  according again   to \e{eq:LTM1B}. 
          \end{pf}

 Let the assumptions of Lemma~\ref{LTM2B}  hold. For the proof of the opposite inclusion $ A \subset A_{0}^{**}$, we have to check relation   \e{eq:LTM1B} for all $ f \in L^2 $ and   $g\in L^2  $  such that ${\cal A} f \in L^2 $ and   ${\cal A}_{*}g\in L^2 $. Suppose now that $v(\xi)\neq 0$ for all $\xi\in{\Bbb R}$.
        Let $\chi_{n}$ be the operator of multiplication by the function $\chi_{n} (x) =\chi(x/n) $ where $\chi=\bar{\chi}$, $\chi(0) =1 $ and the Fourier transform $\hat{\chi}=\Phi\chi\in C_{0}^\infty ({\Bbb R})$. Since $v f \in L^1 $ and   $\chi_{n} \bar{s} g \in L^1 $, it follows from the Fubini theorem that
       \begin{equation}
 ({\cal A} f, \chi_{n}g) =   (f, {\cal A}_{*}\chi_{n} g),\q \forall f,g \in L^2.
\label{eq:LTM3B}\end{equation} 
We  have to pass  here to the limit $n\to \infty$. Since ${\cal A} f\in L^2  $, $g\in L^2  $  and $ \chi_{n} \to I$ strongly in this space, the left-hand side of \e{eq:LTM3B} converges to $({\cal A} f,  g) $.

Let us now consider the right-hand side of \e{eq:LTM3B}. According to \e{eq:SST} we have
\begin{equation}
  {\cal A}_{*}\chi_{n} g=   T_{n}  {\cal A}_{*}  g
\label{eq:LTM3B1}\end{equation}
where
\begin{equation}
T_{n}=  \bar{v}  \Phi \chi_{n} \Phi^* \bar{v}^{-1}.
\label{eq:LTM8B}\end{equation}
Quite similarly to Lemma~\ref{LTM3}, we obtain the following assertion.

\begin{lemma}\label{LTM3B}
Suppose that  $v(\xi)\neq 0$ for all $\xi\in{\Bbb R}$ and that
\begin{equation}
\max_{|\xi-\eta|\leq 1} |v(\eta)|   |v(\xi)|^{-1}   <\infty.                         
\label{eq:LTM8C}\end{equation}
 Then the operators $T_{n}$ defined by formula \e{eq:LTM8B} are bounded in the space $L^2  $, and their norms are bounded uniformly in $n$. Moreover,     $T_{n} \to I$  strongly as $n\to\infty$.
    \end{lemma}
    
    Note that condition \e{eq:LTM8C} admits an exponential decay of the function $v(\xi)$ as $|\xi|\to\infty$, but not a more rapid one. In particular, function \e{eq:LMMv} is allowed.
    
    It follows from equality \e{eq:LTM3B1} and Lemma~\ref{LTM3B} that if $ {\cal A}_{*} g\in L^2$, then $ {\cal A}_{*}\chi_{n} g \to  {\cal A}_{*} g$ as $n\to\infty$. This allows us to pass to the limit $n\to\infty$ in the right-hand side of 
     \e{eq:LTM3B} which yields relation  \e{eq:LTM1B} for all $ f \in L^2 $ and   $g\in L^2  $  such that ${\cal A} f \in L^2 $ and   ${\cal A}_{*}g\in L^2 $. This proves that $ A \subset A_{0}^{**}$.
    
  Let us summarize the results obtained.
     
       \begin{theorem}\label{QFQLB}
  Let  the operator $ {\cal A} $ be defined by formula
\e{eq:BF} on the set $\cal S$. Let the functions $s(x)$ and $v(\xi)$ satisfy assumptions \e{eq:BF1} and \e{eq:BF2} where $k_{0}=0$, respectively.    Then the operator $A_{0} $ in $L^2$ defined on the domain  ${\cal D}(A_{0})=\cal S$ by the equality $A_{0}f= {\cal A}f$ admits the closure.
Suppose additionally that $k_{0}<-1/2$ in \e{eq:BF2} and that the assumptions of Lemma~\ref{LTM3B} are  satisfied.
 Then the closure $\bar{A}_{0}=:A$ is given
   by the same equality $A f={\cal A }f$ on the domain   ${\cal D} (A)$ which  consists of  all $f \in   L^2  $ such that ${\cal A} f\in L^2$.    
  \end{theorem}

%%%%%%%%%%%%%%%%%%%%%%%%%%%%%%%%%%%%%%%%
%%%%%%%%%%%%%%%%%%%%%%%%%%%%%%%%%%%%%%%%

 \end{document}